\documentclass{gtart}

\def\ifplaintex{\expandafter\ifx\csname documentclass\endcsname\relax}

\expandafter\ifx\csname epsfbox\endcsname\relax\input epsf\fi

\ifplaintex 
\else
\fi

\def\gt{{\mathsurround=0pt\it $\cal G\mskip-2mu$eometry \&\ 
$\cal T\!\!$opology}}        

\def\gtp{{\mathsurround=0pt\it $\cal G\mskip-2mu$eometry \&\ 
$\cal T\!\!$opology $\cal P\!$ublications}}  

\def\lognumber#1{\def\thelognumber{#1}}
\def\volumenumber#1{\def\thevolumenumber{#1}}
\def\papernumber#1{\def\thepapernumber{#1}}
\def\volumeyear#1{\def\thevolumeyear{#1}}

\def\pagenumbers#1#2{\def\startpage{#1}\def\finishpage{#2}}
\def\published#1{\def\publishdate{#1}}
\def\proposed#1{\def\theproposer{#1}}
\def\seconded#1{\def\theseconders{#1}}
\def\received#1{\def\receiveddate{#1}}
\def\revised#1{\def\reviseddate{#1}}
\def\accepted#1{\def\accepteddate{#1}}
\def\asciititle#1{\def\theasciititle{#1}}
\def\covertitle#1{\def\thecovertitle{#1}}

\long\def\asciiabstract#1{\long\def\theasciiabstract{#1}}
\def\asciikeywords#1{\def\theasciikeywords{#1}}

\let\\\par\let\thevolumenumber\relax\let\thepapernumber\relax
\let\thevolumeyear\relax\let\thesamplenumber\relax\let\startpage\relax
\let\finishpage\relax\let\publishdate\relax\let\receiveddate\relax
\let\reviseddate\relax\let\accepteddate\relax\let\theasciititle\relax
\let\thecovertitle\relax\let\theasciiauthors\relax
\let\theasciiabstract\relax\let\theasciikeywords\relax
\let\theasciiemail\relax\let\theshortauthors\relax\let\theshorttitle\relax

\long\def\maketitlep{   

\count0=\startpage

\gt\hfill      
\hbox to 77pt{\vbox to 0pt{\vglue -15pt\epsfbox{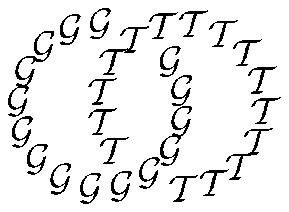}\vss}\hss}
\break
{\small\ifx\thesamplenumber\relax 
Volume \else Sample
\fi\thevolumenumber\ (\thevolumeyear)
\startpage--\finishpage\nl
Published: \publishdate}
\vglue 0.5truein plus 0.4fil minus 0.1truein

{\parskip=0pt\leftskip 0pt plus 1fil\def\\{\par\smallskip}{\ifplaintex\large
\else\Large\fi\bf\thetitle}\par\medskip}   

\vglue 0pt plus 0.1fil 

{\parskip=0pt\leftskip 0pt plus 1fil\def\\{\par}{\sc\theauthors}
\par\medskip}

\vglue 0pt plus 0.1fil 

{\small\parskip=0pt\let\newline\\
{\leftskip 0pt plus 1fil\def\\{\par}{\sl\theaddress}\par}
\expandafter\ifx\theemail\relax    
\relax\else\vglue 5pt plus 0.02fil minus 2pt\def\\{\stdspace{\rm 
and}\stdspace} 
\cl{Email:\stdspace\tt\theemail}\fi
\ifx\theurl\relax                  
\relax\else\vglue 5pt plus 0.02fil minus 2pt\def\\{\stdspace{\rm 
and}\stdspace}
\cl{URL:\stdspace\tt\theurl}\fi\par}

\vglue 7pt plus 0.3fil minus 3pt

{\bf Abstract}
\vglue 5pt plus 0.1fil minus 2pt

\theabstract

\vglue 7pt plus 0.3fil minus 3pt

{\bf AMS Classification numbers}\quad Primary:\quad \theprimaryclass

Secondary:\quad \thesecondaryclass

\vglue 5pt plus 0.3fil minus 2pt

{\bf Keywords:}\quad \thekeywords

\vglue 10pt plus 0.5fil minus 5pt

{\small  Proposed: \theproposer\hfill Received: \receiveddate\nl
Seconded: \theseconders\hfill 
\ifx\reviseddate\relax                         
Accepted: \accepteddate                        
\else
Revised: \reviseddate                          
\fi}
\eject
}       

\let\maketitlepage\maketitlep
\let\maketitle\maketitlepage

\font\phead=cmsl9 scaled 950
\font=cmsl9 scaled 1050
\font\pnum=cmbx10 scaled 913
\font\lnum=cmbx10 
\font\pfoot=cmsl9 scaled 950
\font=cmsl9 scaled 1050
\ifplaintex
\headline{\vbox to 0pt{\vskip -4.5mm\line{\small\phead\ifnum
\count0=\startpage ISSN 1364-0380 (on line)
1465-3060 (printed) \hfill {\pnum\folio}\else\ifodd\count0\def\\{ }%
\ifx\theshorttitle\relax\thetitle\else\theshorttitle\fi\hfill{\pnum\folio}
\else\def\\{ and }{\pnum\folio}\hfill\ifx\theshortauthors\relax\theauthors
\else\theshortauthors\fi\fi\fi}\vss}}
\footline{\vbox to 0pt{\vglue 0mm\line{\small\pfoot\ifnum\count0=\startpage
\copyright\ \gtp\hfill\else
\gt, Volume \thevolumenumber\ (\thevolumeyear)\hfill\fi}\vss
}}
\else
\makeatletter
\def\@oddhead{{\small\lhead\ifnum\count0=\startpage ISSN 1364-0380 (on line)
1465-3060 (printed) \hfill {\lnum\number\count0}\else\ifodd\count0
\def\\{ }\ifx\theshorttitle\relax \thetitle \else\theshorttitle\fi\hfill
{\lnum\number\count0}\else\def\\{ and }{\lnum\number\count0}
\hfill\ifx\theshortauthors\relax 
\theauthors\else\theshortauthors\fi\fi\fi}}\def\@evenhead{@oddhead}
\def\@oddfoot{\small\lfoot\ifnum\count0=\startpage\copyright\ \gtp\hfill\else
\gt, Volume \thevolumenumber\ (\thevolumeyear)\hfill\fi}
\def\@evenfoot{@oddfoot}
\makeatother
\fi

\newwrite\gtoutfile
\long\gdef\makeheadfile{  
{\def\\{, }\def\s{ }
\immediate\openout\gtoutfile head.xxx
\immediate\write\gtoutfile{To: math@arxiv.org}
\immediate\write\gtoutfile{Subject: put OR rep NNNNN:pppp}
\immediate\write\gtoutfile{--text follows this line--}
\immediate\write\gtoutfile{Proxy-for: \ifx\theasciiauthors\relax
\theauthors\else\theasciiauthors\fi\s<\ifx\theasciiemail\relax\theemail\else\theasciiemail\fi>}
\immediate\write\gtoutfile{\noexpand\\}
\immediate\write\gtoutfile{Authors: \ifx\theasciiauthors\relax
\theauthors\else\theasciiauthors\fi}
{\def\\{ }\immediate\write\gtoutfile{Title: \ifx\theasciititle\relax
\thetitle\else\theasciititle\fi}}
\immediate\write\gtoutfile{Subj-class: GT or GR or SG or ...}
\immediate\write\gtoutfile{MSC-class: \theprimaryclass\ifx\thesecondaryclass\relax\else, \thesecondaryclass\fi}
\immediate\write\gtoutfile{Journal-ref: Geom. Topol. \thevolumenumber\s
(\thevolumeyear) \startpage-\finishpage}
\immediate\write\gtoutfile{Comments: Published in Geometry and Topology at}
\immediate\write\gtoutfile{    http://www.maths.warwick.ac.uk/gt/GTVol\thevolumenumber/paper\thepapernumber.abs.html}
\immediate\write\gtoutfile{\noexpand\\}
\immediate\write\gtoutfile{}
\ifx\theasciiabstract\relax
\immediate\write\gtoutfile{\theabstract}\else
\immediate\write\gtoutfile{\theasciiabstract}\fi
\immediate\write\gtoutfile{}
\immediate\write\gtoutfile{\noexpand\\}
\immediate\write\gtoutfile{}
\immediate\closeout\gtoutfile}}  

\def\maketitlepage{\maketitlep\makeheadfile}
\let\maketitle\maketitlepage

\lognumber{219}
\volumenumber{7}\papernumber{2}\volumeyear{2003}
\pagenumbers{33}{90}
\received{15 November 2001}
\revised{6 January 2003}
\accepted{31 Januray 2003}
\published{1 February 2003}
\proposed{Walter Neumann}
\seconded{Benson Farb, David Gabai}

\usepackage[all]{xy}
\usepackage{amsmath}
\usepackage{amssymb}
\input{epsf.tex}

\newcommand{\centeredepsfbox}[1]{\centerline{\epsfbox{#1}}}

\numberwithin{equation}{section}

\newcommand{\nb}[1]{#1--}

\newcommand\ds\displaystyle

\theoremstyle{definition}

\newtheorem*{problem}{Problem}

\theoremstyle{plain}
\newtheorem{theorem}{Theorem}[section]
\newtheorem*{theorem*}{Theorem}
\newtheorem*{corollary*}{Corollary}
\newtheorem{proposition}[theorem]{Proposition}
\newtheorem{lemma}[theorem]{Lemma}
\newtheorem{sublemma}[theorem]{Sublemma}
\newtheorem{corollary}[theorem]{Corollary}
\newtheorem{claim}[theorem]{Claim}

\DeclareMathOperator{\Out}{Out}

\DeclareMathOperator{\Length}{len}

\DeclareMathOperator{\Homeo}{Homeo}

\DeclareMathOperator{\Isom}{Isom}

\DeclareMathOperator\diam{diam}

\DeclareMathOperator\image{image}

\DeclareMathOperator\interior{int}
\DeclareMathOperator\GFL{\mathcal{GFL}}

\DeclareMathOperator\Hull{Hull}

\DeclareMathOperator\Susp{Susp}

\DeclareMathOperator\GF{\mathcal{GF}}
\DeclareMathOperator\domain{domain}

\DeclareMathOperator\genus{genus}

\newcommand\R{{\mathbf R}}
\newcommand\reals{\R}

\newcommand\hyp{\mathbf{H}}

\newcommand\Z{{\mathbf Z}}

\newcommand\solv{{\textsc{solv}}}

\renewcommand\P{\mathcal P}

\renewcommand\H{{\mathcal{H}}}

\newcommand\inject{\hookrightarrow}
\newcommand\homeo{\approx}

\newcommand\Sum{\sum}
\newcommand\infinity{\infty}
\newcommand\bndry{\partial}
\newcommand{\bdy}{\bndry}
\newcommand{\from}{\colon}
\newcommand\composed{\circ}
\newcommand\suchthat{\bigm|}
\newcommand\inverse{{-1}}
\newcommand\inv{\inverse}

\newcommand\union{\cup}

\newcommand\absvalue[1]{\left| #1 \right|}
\newcommand\abs[1]{\absvalue{#1}}
\newcommand\norm[1]{\left\| #1 \right\|}
\newcommand\wt\widetilde

\newcommand\Id{\text{Id}}
\newcommand\A{\mathcal A}
\newcommand\B{\mathcal B}
\newcommand\intersect{\cap}
\newcommand\FP{{\mathcal{FP}}}

\DeclareMathOperator\QD{QD}

\newcommand\restrict{\bigm|}
\newcommand\subgroup{<}
\newcommand\semidirect{\rtimes}

\newcommand\Teichmuller{Teichm\"uller}

\newcommand\MCG{{\mathcal M \mathcal C \mathcal G}}

\newcommand\cross{\times}

\newcommand\F{{\cal F}}
\newcommand\G{{\cal G}}

\newcommand\M{{\cal M}}
\newcommand\Mod\M
\newcommand\ML{{\cal ML}}

\newcommand\PMF{\P\MF}
\newcommand\MF{{\cal MF}}
\newcommand\GL{{\cal GL}}

\renewcommand\S{{\cal S}}

\newcommand\Union\bigcup
\newcommand\C{\mathcal C}

\newcommand\T{{\cal T}}
\newcommand\<\langle
\renewcommand\>\rangle

\newcommand\carries\succ
\newcommand\carriedby\prec
\newcommand\splitsto\carries
\newcommand\collapsesto\carriedby

\newcommand\Fol{f}

\newcommand\ray[2]{\overrightarrow{[#1,#2)}}
\newcommand\geodesic[2]{\overleftrightarrow{(#1,#2)}}

\newcommand\inj{{\mathrm{inj}}}

\begin{document}

\title{Stable \Teichmuller\ quasigeodesics\\and ending laminations}
\covertitle{Stable Teichm\noexpand\"uller quasigeodesics\\and ending laminations}
\asciititle{Stable Teichmueller quasigeodesics and ending laminations}
\author{Lee Mosher}
\address{Deptartment of Mathematics and Computer Science\\Rutgers University,
Newark, NJ 07102}
\email{mosher@andromeda.rutgers.edu}

\begin{abstract} 
We characterize which cobounded quasigeodesics in the \Teichmuller\
space $\T$ of a closed surface are at bounded distance from a
geodesic. More generally, given a cobounded lipschitz path $\gamma$ in
$\T$, we show that $\gamma$ is a quasigeodesic with finite Hausdorff
distance from some geodesic if and only if the canonical hyperbolic
plane bundle over $\gamma$ is a hyperbolic metric space. As an
application, for complete hyperbolic 3--manifolds $N$ with finitely
generated, freely indecomposable fundamental group and with bounded
geometry, we give a new construction of model geometries for the
geometrically infinite ends of $N$, a key step in Minsky's proof of
Thurston's ending lamination conjecture for such manifolds.
\end{abstract}

\asciiabstract{ 
We characterize which cobounded quasigeodesics in the Teichmueller
space T of a closed surface are at bounded distance from a
geodesic. More generally, given a cobounded lipschitz path gamma in
T, we show that gamma is a quasigeodesic with finite Hausdorff
distance from some geodesic if and only if the canonical hyperbolic
plane bundle over gamma is a hyperbolic metric space. As an
application, for complete hyperbolic 3-manifolds N with finitely
generated, freely indecomposable fundamental group and with bounded
geometry, we give a new construction of model geometries for the
geometrically infinite ends of N, a key step in Minsky's proof of
Thurston's ending lamination conjecture for such manifolds.}


\asciikeywords{Teichmueller space, hyperbolic space, quasigeodesics,
ending laminations}

\maketitlepage

\section{Introduction}

\subsection{Stable \Teichmuller\ quasigeodesics}

If $X$ is a geodesic metric space which is hyperbolic in the sense of Gromov, then
quasigeodesics in $X$ are ``stable'': each quasigeodesic line, ray, or segment in $X$ has finite
Hausdorff distance from a geodesic, with Hausdorff distance bounded solely in terms of
quasigeodesic constants and the hyperbolicity constant of $X$. In the case of hyperbolic
$2$--space, stability of quasigeodesics goes back to Morse \cite{Morse:fundamentalclass}. 

Consider a closed, oriented surface $S$ of genus $\ge 2$. Its \Teichmuller\ space
$\T$ is often studied by finding analogies with hyperbolic metric spaces. For example, the
mapping class group $\MCG$ of $S$ acts on $\T$ properly discontinuously by isometries, with
quotient orbifold $\Mod =\T / \MCG$ known as Riemann's moduli space, and $\Mod$ is often viewed
as the analogue of a finite volume, cusped hyperbolic orbifold. Although this analogy is
limited---$\T$ is \emph{not} a hyperbolic metric space in the \Teichmuller\ metric
\cite{MasurWolf:TeichNotHyp}, and indeed there is no $\MCG$ equivariant hyperbolic metric on $\T$
with ``finite volume'' quotient \cite{BrockFarb:curvature}---nonetheless the analogy has recently
been strengthened and put to use in applications \cite{Minsky:quasiprojections},
\cite{MasurMinsky:complex1}, \cite{FarbMosher:quasiconvex}.

Minsky's projection theorem \cite{Minsky:quasiprojections} gives a version in $\T$ of stability
of quasigeodesic segments, assuming coboundedness of the corresponding geodesic.\break  Given a
$\lambda,\eta$ quasigeodesic segment $\gamma$ in $\T$, let $g$ be the geodesic with the same
endpoints as $\gamma$. If $g$ is \emph{cobounded}, meaning that each point of $g$ represents a
hyperbolic structure with injectivity radius bounded away from zero, then according to Minsky's
theorem the Hausdorff distance between $\gamma$ and $g$ is bounded by a constant depending only on
$\lambda$, $\eta$, and~$\epsilon$. On the other hand, Masur and Minsky
\cite{MasurMinsky:unstable} produced a cobounded quasigeodesic line---the orbit of a partially
pseudo-Anosov cyclic subgroup---for which there does not exist any geodesic line at finite
Hausdorff distance.

In joint work with Benson Farb \cite{FarbMosher:quasiconvex} we formulated a theory of convex
cocompact groups of isometries of $\T$, pursuing an analogy with hyperbolic spaces. In the course
of this work we needed to consider the following:
\begin{problem} Give conditions on a cobounded, lipschitz path $\gamma$ in $\T$, in terms of the
geometry of the canonical hyperbolic plane bundle over $\gamma$, which imply that $\gamma$ is a
quasigeodesic with finite Hausdorff distance from a geodesic.
\end{problem} 
Recall that for any path $\gamma \from I \to \T$, $I \subset \reals$ closed and
connected, there is a canonical bundle $\H_\gamma \to I$ of hyperbolic planes over $\gamma$ on
which $\pi_1 S$ acts, so that for each $t \in I$ the marked hyperbolic surface $\S_t = \H_t /
\pi_1 S$ represents the point $\gamma(t)\in\T$. After perturbation of $\gamma$, the fiberwise
hyperbolic metric on $\H_\gamma$ extends to a $\pi_1 S$--equivariant, piecewise Riemannian metric
on $\H_\gamma$, whose large scale geometry is independent of the perturbation; see
\cite{FarbMosher:quasiconvex}, proposition~4.2, recounted below as
proposition~\ref{PropMetricPerturbation}. Thus the above problem asks for an invariant of the
large-scale geometry of $\H_\gamma$ that characterizes when $\gamma$ is a stable quasigeodesic. 

The first goal of this paper is a solution of the above problem:

\begin{theorem}[Hyperbolicity implies stability]
\label{TheoremStableTeichGeod}
If $\gamma \from I \to \T$ is a cobounded, lipschitz map defined on a closed, connected subset $I
\subset \reals$, and if $\H_\gamma$ is a hyperbolic metric space, then $\gamma$ is a
quasigeodesic and $\gamma(I)$ has finite Hausdorff distance from some geodesic $g$ in~$\T$ sharing
the same endpoints as~$\gamma$. To be precise, for every bounded subset $\B \subset \Mod$, and
every $\rho \ge 1$, $\delta\ge 0$ there exists $\lambda\ge 1$, $\eta\ge 0$, $A\ge 0$ such that
the following holds: if $\gamma\from I \to\T$ is $\B$--cobounded and $\rho$--lipschitz, and if
$\H_\gamma$ is $\delta$--hyperbolic, then $\gamma$ is a $\lambda,\eta$ quasigeodesic, and there
exists a geodesic $g \from I \to \T$, with $\gamma(t)=g(t)$ for all $t \in \bdy I$, such that the
Hausdorff distance between $\gamma(I)$ and $g(I)$ is~$\le A$.
\end{theorem} 

\paragraph{Remark} This theorem, and its application to Minsky's results on the ending
lamination conjecture, have been discovered and proved independently by Brian Bowditch
\cite{Bowditch:stacks}, using different methods.

\paragraph{Remark} In the case of a compact segment $I=[a,b]$, the first sentence of the
theorem has no content: $\H_\gamma$ is quasi-isometric to the hyperbolic plane and so is
hyperbolic, and $\gamma$ is Hausdorff equivalent to the \Teichmuller\ geodesic
$\overline{\gamma(a)\gamma(b)}$. The quantifiers in the second sentence are therefore necessary
in order to say anything of substance when $\gamma$ is a segment. Even in the case of a line or
ray, where the first sentence actually has content, the second sentence gives a lot more
information.

\paragraph{Remark} While $\gamma$ is only assumed lipschitz (coarse Lipschitz would do), the
conclusion shows that $\gamma$ is a quasigeodesic. Note that for any lipschitz path $\gamma$ and
geodesic $g$, to say that $\gamma$ is a quasigeodesic at finite Hausdorff distance from $g$ is
equivalent to saying that $\gamma$ and $g$ are asynchronous fellow travellers; see
section~\ref{SectionCoarseGeometry}.
\bigskip

Recent work of Minsky and Rafi provides a strong converse to Theorem~\ref{TheoremStableTeichGeod}
in the case of a bi-infinite path, and putting these together we get the following result:

\begin{corollary}
Given a cobounded, lipschitz path $\gamma \from \reals \to \T$, the following are equivalent:
\begin{enumerate}
\item \label{ItemGeodFellowTraveller}
$\gamma$ is a quasigeodesic and there is a geodesic $g$ at finite Hausdorff distance from
$\gamma$.
\item \label{ItemHypBundle}
$\H_\gamma$ is a hyperbolic metric space.
\item \label{ItemH3Bundle}
$\H_\gamma$ is quasi-isometric to $\hyp^3$.
\end{enumerate}
\end{corollary}

\begin{proof} (\ref{ItemHypBundle})$\implies$(\ref{ItemGeodFellowTraveller}) is
Theorem~\ref{TheoremStableTeichGeod}.
(\ref{ItemGeodFellowTraveller})$\implies$(\ref{ItemH3Bundle}) follows from the Minsky--Rafi
theorem \cite{Minsky:boundedgeometry} which says that $\H^\solv_g$ is quasi-isometric to $\hyp^3$,
together with the fact that $\H_\gamma$ is quasi-isometric to $\H^\solv_g$ (see
Proposition~\ref{PropMetricPerturbation}). (\ref{ItemH3Bundle})$\implies$(\ref{ItemHypBundle}) is
immediate.
\end{proof}

\paragraph{Application: hyperbolic surface group extensions}
Theorem~\ref{TheoremStableTeichGeod} is applied in \cite{FarbMosher:quasiconvex} as follows.
Consider a short exact sequence of finitely generated groups of the form $1 \to
\pi_1(S) \to \Gamma_G \to G \to 1$, determined by a group homomorphism $f \from G \to \MCG(S)
\approx \Out(\pi_1(S))$. One of the main results in \cite{FarbMosher:quasiconvex} says that if
$\Gamma_G$ is a word hyperbolic group, then $G$ maps with finite kernel onto a subgroup of
$\MCG(S)$ acting quasiconvexly on \Teichmuller\ space; Theorem~\ref{TheoremStableTeichGeod}
plays a key role in the proof. 

We remark that when the group $G$ is free then the converse is also
proved in \cite{FarbMosher:quasiconvex}: a free, convex cocompact subgroup $G\subgroup\Isom(\T)$
has a word hyperbolic extension group $\Gamma_G = \pi_1(S) \semidirect G$.

\subsection{Ending laminations}

The second goal of this paper is to apply Theorem~\ref{TheoremStableTeichGeod} to give a new
construction of model manifolds for geometrically infinite ends, used by Minsky in proving
a special case of Thurston's ending lamination conjecture \cite{Minsky:endinglaminations}. 

Let $N$ be a complete hyperbolic \nb{3}manifold satisfying the following properties:
\begin{enumerate}
\item[(1)] The fundamental group $\pi_1 N$ is finitely generated, freely indecomposable.
\item[(2)] The action of $\pi_1 N$ on $\hyp^3$ has no parabolic elements. 
\end{enumerate}
Work of Thurston \cite{Thurston:GeomTop} and Bonahon \cite{Bonahon:ends}
describes the topology and geometry of $N$ nears its ends. The manifold $N$ is the interior of a
compact \nb{3}manifold with boundary denoted $\overline N$, and so there is a one-to-one
correspondence between ends $e$ of $N$ and components $S_e$ of $\bdy\overline N$. Associated to
each end $e$ of $N$ there is an \emph{end invariant}, a geometric/topological structure on the
associated surface $S_e$, which describes the behavior of $N$ in the end~$e$. The end invariant
comes in two flavors. In one case, the end $e$ is \emph{geometrically finite} and the end
invariant is a point in the \Teichmuller\ space of $S_e$. In the case where $e$ is geometrically
infinite, it follows that $e$ is \emph{simply degenerate} and the end invariant is the
\emph{ending lamination}, an element of the geodesic lamination space of $S_e$. 

Thurston's ending lamination conjecture says that $N$ is determined up to isometry by its
topological type and its end invariants: if $N,N'$ are complete hyperbolic \nb{3}manifolds as
above, and if $f\from N \to N'$ is a homeomorphism respecting end invariants, then $f$ is properly
homotopic to an isometry. Minsky's theorem says that this is true if $N,N'$ each have
\emph{bounded geometry}, meaning that injectivity radii in $N$ and in $N'$ are both bounded away
from zero.

The heart of Minsky's result is theorem 5.1 of \cite{Minsky:endinglaminations}, the construction
of a \emph{model manifold} for $N$, a proper geodesic metric space $M$ equipped with a
homeomorphism $f \from N \to M$, such that the metric on $M$ depends only on the end invariants
of $N$, and such that $f$ is properly homotopic to a map for which any lift $\wt N
\to\wt M$ is a quasi-isometry. When $N,N'$ as above have the same end invariants, it follows that
they have isometric model manifolds, and so any homeomorphism $N \to N'$ respecting end
invariants is properly homotopic to a map each of whose lifts $\wt N = \hyp^3 \to \hyp^3 = \wt N'$
is an equivariant quasi-isometry. An application of Sullivan's rigidity theorem
\cite{Sullivan:ErgodicAtInfinity} shows immediately that this quasi-isometry is a bounded distance
from an equivariant isometry, proving the ending lamination conjecture in the presence of
bounded geometry.

Minsky's construction of the model manifold for $N$ is broken into three pieces: a standard
construction in the compact core of $N$; a construction in the geometrically finite ends using
results of Epstein and Marden \cite{EpsteinMarden:convex}; and a construction in the
geometrically infinite ends using Minsky's earlier results \cite{Minsky:GeodesicsAndEnds}.  The
geometrically infinite construction is described in two results. In the singly degenerate case, 
Section~5 of \cite{Minsky:endinglaminations} describes a model manifold for a neighborhood of the
end. In the doubly degenerate case, Corollary~5.10 of \cite{Minsky:endinglaminations}, recalled
below in Theorem~\ref{TheoremELCDD}, describes a model manifold for all of $N$. It is these two
results which we will prove anew, by applying Theorem~\ref{TheoremStableTeichGeod}. 

Given a closed surface $S$, let $\T(S)$ be the \Teichmuller\ space, $\PMF(S)$ the
projective measured foliation space, and $\overline T(S) = T(S) \union \PMF(S)$ Thurston's
compactification. As mentioned above, associated to each geodesic segment, ray, or line $g \from I
\to \T(S)$, where $I$ is a closed connected subset of $\reals$, there is a singular \solv\ metric
on $S\cross I$ denoted $\S^\solv_g$ with universal cover $\H^\solv_g$, such that the induced
conformal structure on $S \cross t$ represents the point $g(t) \in \T(S)$.

Recall that a \emph{pleated surface} in the hyperbolic \nb{3}manifold $N$ is a map denoted $f
\from (S,\sigma)\to N$ where $S$ is a closed surface, $\sigma$ is a hyperbolic structure on $S$,
the map $f$ takes rectifiable paths to rectifiable paths of the same length, and there is a
geodesic lamination $\lambda$ in the hyperbolic structure $\sigma$, called the \emph{pleating
locus} of $f$, such that $f$ is totally geodesic on each leaf of $\lambda$ and on each
complementary component of $\lambda$. Given a $\pi_1$--injective map of a closed surface $S \to
N$, let $\Sigma(S \to N)$ be the set of points $\sigma \in \T(S)$ represented by pleated surfaces
$(S,\sigma) \to N$ in the homotopy class of the map $S \to N$.

Consider now a hyperbolic \nb{3}manifold $N$ as above. An end $e$ of $N$ has neighborhood $N_e
\approx S_e \cross [0,\infinity)$. By Bonahon's theorem, the end $e$ is \emph{geometrically
infinite} if and only if each neighhborhood of $e$ contains the image of a pleated surface in $N$
homotopic to the inclusion $S_e \inject N$. If $\wt N$ is the covering space of $e$ corresponding
to the injection $\pi_1(S_e) \inject \pi_1(N)$, we obtain a homeomorphism $\wt N\approx S_e \cross
(-\infinity,+\infinity)$, with one end corresponding to $N_e$. Assuming $e$ is geometrically
infinite, we say that the inclusion $\pi_1(S_e) \inject \pi_1(N)$ is \emph{singly degenerate} if
$\wt N$ has one geometrically infinite end; otherwise $\wt N$ has two geometrically infinite
ends, in which case we say that $\pi_1(S_e) \inject \pi_1(N)$ is \emph{doubly degenerate}. In the
doubly degenerate case, either the covering map $\wt N \to N$ is degree~1 and so $N \approx S
\cross (-\infinity,+\infinity)$, or $\wt N \to N$ has degree~2 and there is an orbifold fibration
of $N$, with generic fiber $S_e$, whose base space is the \emph{ray orbifold} $[0,\infinity)$ with
a $\Z/2$ mirror group at the point $0$; we'll review these facts in section~\ref{SectionDD}.

\begin{theorem}[Doubly degenerate model manifold]
\label{TheoremELCDD} Let $N$ be a bounded geometry hyperbolic \nb{3}manifold satisfying (1) and
(2) above. Let $e$ be an end, and suppose $\pi_1(S_e) \inject \pi_1(N)$ is doubly degenerate.
There exists a unique cobounded geodesic line $g$ in $\T(S_e)$ such that:
\begin{itemize}
\item[\rm(1)] $\Sigma(S_e \to N)$ is Hausdorff equivalent to $g$ in $\T(S_e)$.
\end{itemize} 
Moreover:
\begin{itemize}
\item[\rm(2)] The homeomorphism $S_e \cross (-\infinity,+\infinity) \approx \wt N$ is properly
homotopic to a map which lifts to a quasi-isometry $\H^\solv_g \to\hyp^3$.
\item[\rm(3)] The ideal endpoints of $g$ in $\PMF(S_e)$ are the respective ending laminations of the
two ends of~$\wt N$. 
\end{itemize}
In the degree~2 case, the order~2 covering transformation group on $\wt N$ acts isometrically on
$\S^\solv_g$ so as to commute with the homeomorphism $\S^\solv_g\approx \wt N$.
\end{theorem}

In the degree~2 case, it follows from the theorem that the degree~2 covering map $\S^\solv_g \to
N$ induces a singular \solv\ metric on $N$.

In the doubly degenerate case we have constructed the entire model manifold, and therefore by
applying Sullivan's rigidity theorem \cite{Sullivan:ErgodicAtInfinity} we obtain a complete proof
of the ending lamination conjecture for the case $N\approx S \cross (-\infinity,+\infinity)$ with
two geometrically infinite ends and with bounded geometry. 

In the singly degenerate case, one still needs the arguments of Minsky to construct the model
manifold out of its pieces, before applying Sullivan's rigidity theorem. Here is the theorem
describing the piece of the model manifold corresponding to a singly degenerate end.

\begin{theorem}[Singly degenerate model manifold]
\label{TheoremELC}
Let $N$ be a hyperbolic \nb{3}manifold of bounded geometry satisfying (1) and (2) above, and let
$e$ be a geometrically infinite end. There exists a cobounded geodesic ray $g$ in $\T(S_e)$,
unique up to choice of its finite endpoint, such that:
\begin{itemize}
\item[\rm(1)] $\Sigma(S_e \to N)$ and $g$ are Hausdorff equivalent in $\T(S_e)$.
\end{itemize}
Moreover, $g$ satisfies the following:
\begin{itemize}
\item[\rm(2)] The homeomorphism $S^\solv_g \approx S_e \cross [0,\infinity)\approx N_e$ is properly
homotopic to a map which lifts to a quasi-isometry of universal covers $\H^\solv_g\to\wt N_e$.
\item[\rm(3)] The ideal endpoint of $g$ in $\PMF(S_e)$ is the ending lamination of $e$.
\end{itemize}
\end{theorem}

Minsky's proofs of these results in \cite{Minsky:endinglaminations} depend on pleated surface
arguments found in section~4 of that paper, as well as results found in
\cite{Minsky:GeodesicsAndEnds} concerning projection arguments in \Teichmuller\ space and
harmonic maps from Riemann surfaces into hyperbolic \nb{3}manifolds. In particular, the
statements (1) in these two theorems are parallel to harmonic versions stated in Theorem~A of
\cite{Minsky:GeodesicsAndEnds}.

Our proofs are found in section~\ref{SectionELC}; here is a sketch in the
doubly degenerate case $N \homeo S \cross (-\infinity,+\infinity)$. The first part uses pleated
surface arguments, several borrowed from \cite{Minsky:endinglaminations} section~4. There is a
sequence of pleated surfaces $f_i \from (S,\sigma_i)\to N$, ($i=-\infinity,\ldots,+\infinity$),
whose images progress in a controlled fashion from $-\infinity$ to $+\infinity$ in $N$. This
allows one to construct a cobounded, lipschitz path $\gamma \from \reals \to \T(S)$ with
$\gamma(i)=\sigma_i$, and to
construct a map $\S_\gamma \to N$ in the correct homotopy class each of whose lifts to the
universal covers $\H_\gamma \to \hyp^3$ is a quasi-isometry. 

For the next part of the proof, since $\hyp^3$ is Gromov hyperbolic, so is $\H_\gamma$, and so
Theorem~\ref{TheoremStableTeichGeod} applies to the path $\gamma$. The output is a cobounded
\Teichmuller\ geodesic $g$ such that $\gamma$ and $g$ are asynchronous fellow travellers. It is
then straightforward to show that $g$ satisfies the desired conclusions of
Theorem~\ref{TheoremELCDD}.

The singly degenerate case is similar, requiring an additional argument to show that $N_e$
itself has Gromov hyperbolic universal cover, so that Theorem~\ref{TheoremStableTeichGeod} can
apply to produce the desired \Teichmuller\ ray. 

\bigskip

Finally, suppose that $N$ satisfies (1) and (2) but we do not assume that $N$ has bounded
geometry. Suppose however that some geometrically infinite end $e$ of $N$ does have bounded
geometry, in the sense that for some $\epsilon>0$ each point of $N_e$ has injectivity radius $\ge
\epsilon$. In this case we are still able to construct a model geometry for $N_e$: there is a
geodesic ray $g \from [0,\infinity) \to \T$ and a map $\S_g \to N_e$ in the correct homotopy
class which lifts to a quasi-isometry $\H_g \to \wt N_e$, and the other conclusions of
Theorem~\ref{TheoremELC} hold as well. This case is not covered by Minsky's results in
\cite{Minsky:endinglaminations}, although probably with some work the original construction could
be pushed through. Our construction in this new case is given in
Section~\ref{SectionEndsBoundedGeom}. Thanks to Jeff Brock for asking this question.

\subsection{An outline of the proof of Theorem \ref{TheoremStableTeichGeod}}

We shall outline the proof for a path $\gamma$ whose domain is the whole real line,
$\gamma\from\reals\to \T$. As remarked earlier, $\pi_1 S$ acts on each fiber $\H_t$ of the
hyperbolic plane bundle $\H_\gamma \to \reals$, with quotient a marked hyperbolic surface $\S_t =
\H_t / \pi_1 S$. These surfaces fit together to form the canonical marked hyperbolic surface
bundle $\S_\gamma \to\reals$ over $\gamma$. Note that $\H_\gamma$ is the universal cover of
$\S_\gamma$, with deck transformation group $\pi_1 S$.

The proof uses the Bestvina--Feighn flaring
condition for $\H_\gamma$, a necessary and sufficient condition for Gromov hyperbolicity of
$\H_\gamma$.  The key idea of the proof is to use hyperbolicity of $\H_\gamma$, via the
Bestvina--Feighn flaring condition, to construct ending laminations for $\S_\gamma$, one
lamination $\lambda^-$ for the negative end and another $\lambda^+$ for the positive end. These
two laminations can then be used to construct the desired \Teichmuller\ geodesic~$g$.

The flaring condition is concerned with quasihorizontal paths in $\H_\gamma$, that is, sections
$\ell\from \reals \to \H_\gamma$ of the projection map $\H_\gamma \to \reals$ such that
$\ell$ satisfies a coarse Lipschitz condition. The flaring condition says, roughly, that for any
two quasihorizontal paths $\ell,\ell'\from\reals \to\H_\gamma$, the sequence of distances $d_i =
d_{\H_i}(\ell(i),\ell'(i))$ satisfies an exponential growth property in either forward or
backward time, possibly both.

The central concept in the proof is an exponential growth
property for measured geodesic laminations. Each element $\lambda$ in $\MF$, the space of
measured foliations on $S$, can be represented as a measured geodesic lamination $\lambda_s$ in
each hyperbolic surface $\S_s$, $s \in\reals$. Although not strictly necessary for the proof, it
is nevertheless convenient to arrange that as $s$ varies in $\reals$ the lamination
$\lambda_s$ varies nicely in the bundle $\S$, filling out a 2--dimensional lamination on $\S$; the
proof of this fact relies on basic tools of partially hyperbolic dynamics to show that the
geodesic flows on the surfaces $\S_s$ can be packaged together in a nice manner. The Thurston
length of $\lambda_s$ is obtained by integrating the transverse measure on $\lambda_s$ against the
Lebesgue measure along leaves of $\lambda_s$, yielding a number $\ell_s > 0$. We are interested in
the growth properties of the sequence $\ell_i$, $i \in \Z$. By applying the Bestvina--Feighn
flaring property for $\H_\gamma$, together with an argument using Fubini's theorem, we show that
$\ell_i$ satisfies an exponential growth property in either forward or backward time, possibly
both. As a consequence, each element $\lambda \in \MF$ satisfies the following trichotomy:
\begin{description}
\item[$\lambda$ is realized at $-\infinity$,] meaning that $\ell_i$ goes exponentially to zero as
$i \to -\infinity$ and exponentially to infinity as $i \to +\infinity$; \emph{or}
\item[$\lambda$ is realized at $+\infinity$,] meaning that $\ell_i$ goes exponentially to zero as
$i \to +\infinity$ and exponentially to infinity as $i \to -\infinity$; \emph{or} 
\item[$\lambda$ is finitely realized,] meaning that $\ell_i$ goes exponentially to infinity as $i
\to -\infinity$ or $+\infinity$, and $\ell_i$ is minimized on a subinterval of $\reals$ of
uniformly bounded length.
\end{description}
We also prove that the position at which $\lambda$ is realized, either at $-\infinity$, or at
$+\infinity$, or on a certain subinterval of $\reals$ of uniformly bounded length, is a coarsely
continuous function of $\lambda \in \MF$. By considering the shortest geodesic in each $\S_i$,
letting $i \to \pm\infinity$, and passing to limits in the projective measured foliation space
$\PMF$, we exhibit the existence of geodesic laminations $\lambda^-$, $\lambda^+$ realized at
$-\infinity$, $+\infinity$ respectively. These are the \emph{ending laminations} of the
hyperbolic surface bundle $\S_\gamma$. 

The next step is to prove that $\lambda^-,\lambda^+$ each fill the surface $S$. This uses a trick
that I remember learning in Thurston's class at Princeton, sometime from 1979 to 1983. If, say,
$\lambda^+$ does not fill, then there is a simple closed geodesic $c$ that, together with certain
boundary leaves of $\lambda^+$, cobounds a ``crown surface'' (see Figure~\ref{FigureCrown}). But
then one can play off the exponential decay of $\lambda^+$ against the exponential growth of $c$
to get a contradiction.

It is easy to show that $\lambda^-,\lambda^+$ are topologically inequivalent, and it follows that
they jointly fill the surface, which is precisely the condition needed to exhibit a \Teichmuller\
geodesic $g$ with horizontal and vertical measured foliations equivalent to
$\lambda^-,\lambda^+$. 

It remains to show that $\gamma$ and $g$ are at finite Hausdorff distance. Our first proof of
this was somewhat laborious, requiring one to go through the whole construction of $g$ and check
various additional properties along the way. But then we discovered a compactness argument that
establishes finite Hausdorff distance seemingly by magic. The idea is to fix a compact subset $\B
\subset \Mod$ and constants $\rho,\delta \ge 0$, and to consider the space
$\Gamma_{\B,\rho,\delta}$ consisting of all triples $(\gamma,\lambda^-,\lambda^+)$ such that
$\gamma \from \reals \to \T$ is a $\B$--cobounded, $\rho$--lipschitz path for which
$\H_\gamma$ is $\delta$--hyperbolic, and $\lambda^-,\lambda^+$ are ending laminations for
$\S_\gamma$, normalized to have Thurston length~1 in $\S_0$. The
mapping class group $\MCG$ acts on $\Gamma_{\B,\rho,\delta}$, and we prove that this action is
cocompact, using compact open topology for $\gamma$ and the topology on $\MF$ for
$\lambda^-,\lambda^+$. The proof of cocompactness uses the Ascoli--Arzela theorem together with
the fact that $\delta$--hyperbolicity is closed in the Gromov--Hausdorff topology on metric spaces.
We also prove that if $g$ is the \Teichmuller\ geodesic determined by
$\lambda^-,\lambda^+$, then the point $g(0) \in \T$ is a continuous function of the data
$(\gamma,\lambda^-,\lambda^+)$. It follows that the distance in \Teichmuller\ space between
$\gamma(0)$ and $g(0)$ is a continuous function of the data, and since $\MCG$ acts cocompactly on
the data, then $d_\T(\gamma(0),g(0))$ is uniformly bounded. Suppressed in this exposition is a
delicate parameterization issue, which comes up in proving that $\gamma$ is a quasigeodesic: we
actually prove not just that $\gamma,g$ are at finite Hausdorff distance, but that $\gamma,g$ have
quasigeodesic parameterizations so that $d_\T(\gamma(t),g(t))$ is uniformly bounded.

Some additional work is needed to prove the theorem when $\gamma$ is parameterized by a ray or a
finite interval; in the exposition, these cases are smoothly integrated with the case of a line.

\paragraph{Acknowledgements} Thanks to Benson Farb, for a rich collaboration from which many
offshoots can grow, and in particular for many shared ideas which contributed to this paper.

Thanks also to the Technion in Haifa, Israel, and to Michah Sageev, for hospitality and support
during the June 2000 Conference on Geometric and Combinatorial Group Theory, where this work was
carried out and was first presented. 

Thanks to Amie Wilkinson for telling me how to prove Lemma~\ref{LemmaGFConnection}, and to Jeff
Brock for suggesting the result in Section~\ref{SectionEndsBoundedGeom}.

For useful suggestions in preparing the revised version, thanks to the referee, and to Hossein
Namazi and other members of the Complex Analysis and Geometry Seminar at SUNY Stony Brook. 

This work was partially supported by grants from the National Science Foundation.

\section{Preliminaries}

\subsection{Coarse geometry}
\label{SectionCoarseGeometry}

Consider two metric spaces $X,Y$ and a map $f \from X \to Y$. The map $f$ is \emph{$K,C$ coarse
lipschitz} with $K \ge 1, C \ge 0$ if 
$$d_Y(f(x),f(x')) \le K d_X(x,x') + C \quad\text{for}\quad x,x' \in X.
$$ 
If $f$ is coarse lipschitz, then we say in addition that $f$ is
\emph{uniformly proper} with respect to a proper, monotonic function $\rho \from [0,\infinity)
\to [0,\infinity)$ if 
$$d_Y(f(x),f(x')) \ge \rho(d_X(x,x')) \quad\text{for}\quad x,x' \in X.
$$
The function $\rho$ is called a \emph{properness gauge} for $f$. We say that $f$ is a \emph{$K,C$
quasi-isometric embedding} if $f$ is $K,C$ coarse lipschitz and uniformly proper with
properness gauge $\rho(d) = \frac{1}{K} d - C$. The map $f$ is \emph{$C$--coarsely
surjective} if for all $y\in Y$ there exists $x \in X$ such that $d_Y(f(x),y) \le C$. The map $f$
is a
\emph{$K,C$ quasi-isometry} if it is a $C$--coarsely surjective, $K,C$ quasi-isometric embedding.

A metric space $X$ is \emph{geodesic} if for any $x,y$ there is a rectifiable path $p$ from $x$
to $y$ such that $\Length(p)=d(x,y)$. $X$ is \emph{proper} if closed balls are compact.

A map $\bar f \from Y \to X$ is a \emph{$C$--coarse inverse} for $f$ if $d_{\sup}(\bar f \composed
f,\Id_X)  \le C$ and $d_{\sup}(f\composed \bar f,\Id_Y) \le C$.

Here are some basic facts concerning these concepts:
\begin{lemma} 
\label{LemmaQIFacts}
\quad
\begin{enumerate}
\item 
A coarse lipschitz map is a quasi-isometry if and only if it has a coarse inverse
which is coarse lipschitz.
\item \label{ItemUnifPropIsQI}
Suppose $X,Y$ are geodesic metric spaces. Any coarsely surjective, uniformly proper
map $f \from X \to Y$ is a quasi-isometry. \qed
\end{enumerate}
\end{lemma}
In each of these facts, the constants implicit in the conclusion of the statement depend only
on the constants in the hypothesis.

\paragraph{Quasigeodesics} Given a geodesic metric space $X$, a
\emph{$\lambda,\epsilon$ quasigeodesic} in $X$ is a $\lambda,\epsilon$ quasi-isometric embedding
$\gamma \from I \to X$, where $I$ is a closed, connected subset of $X$. With $I$ is a compact
interval we have a \emph{quasigeodesic segment}, when $I$ is a half-line we have a
\emph{quasigeodesic ray}, and when $I=\reals$ we have a \emph{quasigeodesic line}. 

Two paths $\gamma \from I \to X$, $\gamma' \from I' \to X$ are \emph{asynchronous fellow
travellers} with respect to a $K,C$ quasi-isometry $\phi \from I \to I'$ if there is a constant
$A$ such that $d(\gamma'(\phi(t)),\gamma(t)) \le A$ for $t \in I$. 

Recall that the Hausdorff distance between two sets $A,B \subset X$ is the
infimum of $r \subset \R_+ \union +\infinity$ such that $A$ is contained in the $r$--neighborhood
of $B$, and $B$ is contained in the $r$--neighborhood of $A$.

Consider two paths $\gamma \from I \to X$, $\gamma'\from I' \to X$ such that $\gamma$ is a
quasigeodesic. In this situation, the paths $\gamma,\gamma'$ are asynchronous fellow travellers
if and only if $\gamma'$ is a quasigeodesic and the sets $\gamma(I)$, $\gamma'(I')$ have finite
Hausdorff distance in $X$; moreover, the constants implicit in these properties are uniformly
related. To be precise, suppose that $\gamma$ is a $\lambda,\epsilon$ quasigeodesic. If
$\gamma,\gamma'$ are asynchronous fellow travellers with constants $K,C,A$ as above then
$\gamma'$ is a $\lambda',\epsilon'$ quasigeodesic with $\lambda',\epsilon'$ depending only on
$\lambda,\epsilon,K,C,A$. Conversely, if $\gamma'$ is a $\lambda',\epsilon'$ quasigeodesic then
there exist constants $K,C,A$ depending only on $\lambda,\epsilon,\lambda',\epsilon'$, such that
if $\gamma(I)$, $\gamma'(I')$ have Hausdorff distance $\le A$ then any map $\phi \from I \to I'$
with the property that $d(\gamma(t),\gamma'(\phi(t)) \le A$ is a $K,C$--quasi-isometry.

\subsection{Surface geometry and topology}

Fix a closed, oriented surface $S$ of genus $\ge 2$. 

We review \Teichmuller\ space and the accompanying structures: the mapping class group; measured
foliations and the Thurston boundary; measured geodesic laminations; quadratic differentials,
geodesics, and the \Teichmuller\ metric; and canonical bundles over \Teichmuller\ space. Much of
the material in this section is covered in more detail in Sections~2 and~4 of
\cite{FarbMosher:quasiconvex}.

\paragraph{\Teichmuller\ space and mapping class group}
Let $\Homeo(S)$ be the group of homeomorphisms of $S$ and let $\Homeo_0(S)$ be the normal subgroup
of homeomorphisms isotopic to the identity. The \emph{mapping class group} of $S$ is $\MCG =
\Homeo(S) / \Homeo_0(S)$. Let $\C$ be the set of essential simple closed curves on $S$ modulo
isotopy, that is, modulo the action of $\Homeo_0(S)$. The \emph{\Teichmuller\ space} of $S$,
denoted $\T$, is the set of hyperbolic structures on $S$ modulo isotopy, or equivalently the set
of conformal structures modulo isotopy. There are natural actions of $\MCG$ on $\C$ and on $\T$.
The length pairing $\T\cross\C \to \reals_+$, associating to each $\sigma\in\T$, $C \in \C$ the
length of the unique closed geodesic on the hyperbolic surface $\sigma$ in the isotopy class $C$,
induces an $\MCG$--equivariant embedding $\T\to [0,\infinity)^\C$, giving $\T$ the
$\MCG$--equivariant structure of a smooth manifold of dimension $6g-6$ diffeomorphic to
$\reals^{6g-6}$. The action of $\MCG$ on $\T$ is properly discontinuous and noncocompact, and so
the \emph{moduli space} $\Mod = \T / \MCG$ is a smooth, noncompact orbifold of dimension $6g-6$.
The action of $\MCG$ on $\T$ is faithful except in genus~2 where there is a $\Z/2$ kernel
generated by the hyperelliptic involution. 

\paragraph{Measured foliations and Thurston's boundary}
A \emph{measured foliation} on $S$ is a foliation with finitely many singularities, equipped with
a positive transverse Borel measure, such that each singularity is an \emph{$n$--pronged
singularity} for some $n \ge 3$, modelled on the singularity at the origin of the horizontal
foliation of the quadratic differential $z^{n-2} dz^2$. Given a measured foliation, collapsing a
saddle connection---a leaf segment connecting two singularities---results in another measured
foliation. The \emph{measured foliation space} of $S$, denoted $\MF$, is the set of measured
foliations modulo the equivalence relation generated by saddle collapses
and isotopies. There is a natural action of $\MCG$ on $\MF$. The \emph{geometric intersection
number} pairing $\MF \cross \C \to [0,\infinity)$ assigns to each $\F \in \MF$, $C \in \C$ the
number 
$$\< \F, C\> = \inf_{c \in C, f \in \F} \int_{c} f
$$
where $\int_c f$ is the integral of the transverse measure on $f$ pulled back to a measure on $c$.
This pairing induces $\MCG$--equivariant embedding $i \from \MF \inject [0,\infinity)^\C$.
Multiplying transverse measures by positive real numbers defines an action of $(0,\infinity)$ on
$\MF$, whose orbit space is defined to be $\PMF$. The embedding $i \from
\MF \to [0,\infinity)^\C$ induces an embedding $\PMF \to \P[0,\infinity)^\C$, whose image is
homeomorphic to a sphere of dimension $6g-5$. The composed map $\T \to [0,\infinity)^\C \to
\P[0,\infinity)^\C$ is an embedding, the closure of whose image is a closed ball of dimension
$6g-6$ with interior $\T$ and boundary sphere $\PMF$, called the \emph{Thurston compactification}
$\overline\T =\T\union\PMF$. 

Given a simple closed curve $c$ on $S$, by thickening $c$ to form a foliated annulus with total
transverse measure~$1$, and then collapsing complementary components of the annulus to a spine,
we obtain a measured foliation on $S$ well-defined in $\MF$. This gives an embedding $\C \inject
\MF$, whose induced map $\reals_+ \cross \C \to \MF$ is also an embedding with dense image. With
respect to this embedding, the geometric intersection number function
$\MF\cross(\reals_+\cross\C)\to [0,\infinity)$ extends continuously to an intersection number
$\MF\cross\MF \to[0,\infinity)$, denoted $\<\F_1,\F_2\>$, $\F_1,\F_2 \in \MF$.

A pair of measured foliations is \emph{transverse} if they have the same singular set, they are
transverse in the usual sense at each nonsingular point, and for each singularity $s$ there
exists $n \ge 3$ such that the two foliations are modelled on the horizontal and vertical
foliations of $z^{n-2} dz^2$. Given two points $\F_1,\F_2 \in \MF$, we say that $\F_1,\F_2$
\emph{jointly fill} if for each $\G\in \MF$ we have $\<\F_1,\G\> \ne 0$ or $\<\F_2,\G\> \ne 0$. A
pair $\F_1,\F_2$ jointly fills if and only if they are represented, respectively, by a transverse
pair $f_1,f_2$; moreover, the pair $f_1,f_2$ is unique up to joint isotopy, meaning that if
$f'_1,f'_2$ is any other transverse pair representing $\F_1,\F_2$ then there exists $h
\in\Homeo_0(S)$ such that $f'_1=h(f_1)$ and $f'_2=h(f_2)$. 

The set of jointly filling pairs forms an open subset $\FP \subset \MF \cross \MF$. The image of
this set in $\PMF\cross\PMF$ we denote $\P\FP$.

\paragraph{Measured geodesic laminations}
For details of measured geodesic laminations see \cite{CassonBleiler}, \cite{Thurston:GeomTop}.
Here is a brief review.

Given a hyperbolic structure $\sigma$ on $S$, a \emph{geodesic lamination} on $\sigma$ is a
closed subset of $S$ decomposed into complete geodesics of~$\sigma$. A
\emph{measured geodesic lamination} is a geodesic lamination equipped with a positive, transverse
Borel measure. The set of all measured geodesic laminations on $\sigma$ is denoted
$\ML(\sigma)$.  A measured geodesic lamination, when lifted to the universal cover $\wt\sigma
\approx\hyp^2$ with boundary circle $S^1$, determines a positive Borel measure on the complement
of the diagonal in $S^1 \cross S^1$. This embeds $\ML(\sigma)$ into the space of positive Borel
measures, allowing us to impose a topology on $\ML(\sigma)$ using the weak$^*$ topology on
measures. For another description of the same topology, given any simple closed curve $c$ on $S$
and measured geodesic lamination $\mu$ we define the intersection number $\<c,\mu\>$ by pulling
the transverse measure on $\mu$ back to the domain of $c$ and integrating, and we obtain an
embedding $\ML(\sigma)\to [0,\infinity)^\C$ which is a homeomorphism onto its image, using the
product topology on
$[0,\infinity)^\C$.

The space $\ML(\sigma)$ depends naturally on the hyperbolic structure $\sigma$ in the following
sense. For any two hyperbolic structures $\sigma,\sigma'$ we have a homeomorphism $\ML(\sigma)
\to \ML(\sigma')$, obtained by using the natural identification of the circles at infinity of the
universal covers $\wt\sigma$ and $\wt\sigma'$, via the Gromov boundary of the group
$\pi_1 S$. To visualize this homeomorphism, if $\mu$ is a measured geodesic lamination on $\sigma$
then we may regard $\mu$ as a measured \emph{non}-geodesic lamination on $\sigma'$, which may be
straightened to form the corresponding measured geodesic lamination on~$\sigma'$. When
$\sigma=\sigma'$ this homeomorphism is the identity; and the composition
$\ML(\sigma) \to \ML(\sigma') \to \ML(\sigma'')$ agrees with the map $\ML(\sigma) \to
\ML(\sigma'')$ for three hyperbolic structures $\sigma,\sigma',\sigma''$. We may therefore
identify all of the spaces $\ML(\sigma)$ to a single space denoted~$\ML$.

There is a natural isomorphism $\ML \approx \MF$, obtained by taking a measured
geodesic lamination $\lambda$, collapsing the components of $S-\lambda$ to get a foliation
$\Fol$, and pushing the transverse measure on $\lambda$ forward under the collapse map to get a
transverse measure on $\Fol$; the inverse map takes a measured foliation $\Fol$ and straightens
its leaves to get a geodesic lamination $\lambda$ which collapses back to $\Fol$, and the
transverse measure on $\Fol$ is pulled back under the collapse to define the transverse measure
on $\lambda$. Under this isomorphism, the embeddings of $\MF$ and $\ML$ into
$[0,\infinity)^\C$ agree. 

Since each element of $\MF$ is uniquely represented by a measured geodesic lamination
once the hyperbolic structure is chosen, we will often use the same notation to represent
either an element of $\MF$ or its representative measured geodesic lamination, when the
hyperbolic structure is clear from the context.

Given a hyperbolic structure $\sigma$ on $S$ and a measured geodesic lamination $\mu$ on
$\sigma$, let $d\mu^\perp$ denote the transverse measure on $\mu$, let $d\mu^\parallel$ denote
the leafwise Lebesgue measure on leaves, and let $d\mu = d\mu^\perp \cross d\mu^\parallel$
denote the measure on $S$ obtained locally as the Fubini product of $d\mu^\perp$ with
$d\mu^\parallel$. The support of $d\mu$ is
$\mu$, and the \emph{length of $\mu$ with respect to $\sigma$} is defined 
to be
$$\Length_\sigma(\mu) = \int d\mu.
$$
We need the well known fact that length defines a continuous function: 
\begin{align*}
\T\cross\MF &\mapsto (0,\infinity) \\
 (\sigma,\mu) &\to \Length_\sigma(\mu)
\end{align*}
Consider now a hyperbolic structure $\sigma$ on $S$, a pair $\F_1,\F_2 \in \MF$, and measured
geodesic laminations $\lambda_1,\lambda_2$ on $\sigma$ representing $\F_1,\F_2$ respectively. The
intersection number $\<\F_1,\F_2\>$ has the following interpretation. There exist unique maximal
closed sublaminations $\lambda'_1 \subset\lambda_1,\lambda'_2 \subset\lambda_2$, possibly empty,
with the property that $\lambda'_i$ is transverse to $\lambda_j$, $i \ne j \in \{1,2\}$. By taking
the Fubini product of the transverse measures on $\lambda'_1,\lambda'_2$ and integrating over $S$,
we obtain $\<\F_1,\F_2\>$. Joint filling also has an interpretation: the pair $(\F_1,\F_2)$
jointly fills if and only if $\lambda_1,\lambda_2$ are transverse and each component of
$S-(\lambda_1\union\lambda_2)$ is simply connected.

Associated to a hyperbolic structure $\sigma$ on $S$ we also have the space $\GL(\sigma)$ of
(unmeasured) geodesic laminations with the Hausdorff topology. These spaces also depend naturally
on $\sigma$, and hence we may identify them to obtain a single space $\GL$ depending only on $S$.

\paragraph{Quadratic differentials} Given a conformal structure on $S$, a \emph{quadratic
differential} associates to each conformal coordinate $z$ an expression $q(z) dz^2$ with $f$
holomorphic, such that whenever $z,w$ are two overlapping conformal coordinates we have $q(z) =
q(w) \left(\frac{dw}{dz}\right)^2$. The \emph{area form} of $q$ is expressed in a conformal
coordinate $z=x+iy$ as $\abs{q(z)} \abs{dx} \abs{dy}$, and the integral of this form is a
positive number $\norm{q}$ called the \emph{area}. We say that $q$ is \emph{normalized} if
$\norm{q}=1$.

Given a conformal structure $\sigma$, the Riemann--Roch theorem says that the 
quad\-ratic
differentials on $\sigma$ form a vector space $\QD_\sigma$ of dimension $6g-6$, and as $\sigma$
varies over $\T$ these vector spaces fit together to form a vector bundle $\QD\to\T$. The
normalized quadratic differentials form a sphere bundle $\QD^1 \to \T$. 

Given a quadratic differential $q$, for each point $p \in S$, there exists a conformal coordinate
$z$ in which $p$ corresponds to the origin, and a unique $n \ge 2$, such that $q(z)=z^{n-2}$; if
$n=2$ then $p$ is a \emph{regular point} of $q$, and otherwise $p$ is a zero of order $n-2$. The
coordinate $z$ is unique up to multiplication by $n^{\text{th}}$ roots of unity, and is called a
\emph{canonical coordinate} at~$p$.

Associated to each quadratic differential $q$ is a transverse pair of measured foliations, the
\emph{horizontal measured foliation} $\Fol_x(q)$ and the \emph{vertical measured foliation}
$\Fol_y(q)$, characterized by the property that for each regular canonical coordinate $z=x+iy$,
$\Fol_x(q)$ is the foliation by lines parallel to the $x$--axis with transverse measure $\abs{dy}$,
and $\Fol_y(q)$ is the foliation by lines parallel to the $y$--axis with transverse measure
$\abs{dx}$. As mentioned earlier, for the quadratic differential $z^{n-2} dz^2$ on the complex
plane the horizontal and vertical measured foliations have singularities at the origin which are
models for a transverse pair of $n$--pronged singularities.

Conversely, if $\Fol_x,\Fol_y$ is a transverse pair of measured foliations then there exists a
unique conformal structure $\sigma(\Fol_x,\Fol_y)$ and quadratic differential
$q=q(\Fol_x,\Fol_y)$ such that $\Fol_x=\Fol_x(q)$ and $\Fol_y=\Fol_y(q)$. In particular,
associated to each jointly filling pair $\F_1,\F_2 \in \MF\cross\MF$ are uniquely defined points
$\sigma=\sigma(\F_1,\F_2) \in \T$, $q(\F_1,\F_2) \in \QD_\sigma$. 

We obtain an injective map $\QD\to\MF\cross\MF$ with image $\FP$, given by $q \mapsto
\bigl([\Fol_x(q)],[\Fol_y(q)]\bigr)$, and this map is a homeomorphism between $\QD$
and~$\FP$~\cite{HubbardMasur:qd}.

\paragraph{Geodesics and metric on $\T$}
Associated to each quadratic differential $q$ is a \emph{\Teichmuller\ geodesic} $g \from \reals
\to \T$ defined as follows:
$$g(t) = \sigma(e^{-t} \Fol_x(q),e^{t} \Fol_y(q)), \quad t \in \reals
$$
\Teichmuller's theorem says that any two points $p \ne q \in \T$ lie on a \Teichmuller\ line $g$,
and that line is unique up to an isometry of the parameter line $\reals$. Moreover, if
$p=g(s)$ and $q=g(t)$, then the formula $d(p,q) = \abs{s-t}$ defines a proper, geodesic metric on
$\T$, called the \emph{\Teichmuller\ metric}. The \emph{positive ending lamination} of the
\Teichmuller\ geodesic $g$ is
defined to be the point $\P\Fol_y(q)\in\PMF$, and the \emph{negative ending lamination} of $g$ is
$\P\Fol_x(q)$. By uniqueness of representing transverse pairs as described above, it follows that
$\image(g)$ is completely determined by the pair of points $\P\xi=\P\Fol_x(q), \P\eta=\P\Fol_y(q)
\in \PMF$, and we write
$$\image(g) = \geodesic{\P\xi}{\P\eta}.
$$
Given any $\sigma \in \T$ and $\P\eta \in \PMF$, there is a unique geodesic ray denoted
$\ray{\sigma}{\P\eta}$ with finite endpoint $\sigma$, given by the formula above with $t \ge 0$.

The group $\MCG$ acts isometrically on $\T$, and so the metric on $\T$ descends to a
proper geodesic metric on the moduli space $\Mod$. A subset $A \subset \T$ is said to be
\emph{cobounded} if there exists a bounded subset $B \subset \T$ whose translates under $\MCG$
cover $A$; equivalently, the projection of $A$ to $\Mod$ has bounded image. In a similar way we
define \emph{cocompact} subsets of~$\T$. A subset of $\T$ is cocompact if and only if it is
closed and cobounded.  Usually we express coboundedness in terms of some bounded subset $\B
\subset \Mod$; a subset $A\subset \T$ is said to be \emph{$\B$--cobounded} if the projection of
$A$ to $\Mod$ is contained in $\B$.

Mumford's theorem provides a gauge for coboundedness. Given $\epsilon>0$, define
$\T_\epsilon$ to be the set of hyperbolic structures $\sigma$ whose shortest closed geodesic has
length $\ge\epsilon$, and define $\Mod_\epsilon$ to be the projected image of $\T_\epsilon$.
Mumford's theorem says that the sets $\Mod_\epsilon$ are all compact, and their
union as $\epsilon\to 0$ is evidently all of $\Mod$. It follows that a subset $A \subset \T$ is
cobounded if and only if it is contained in some $\T_\epsilon$. We generally will not rely
Mumford's gauge, instead relying on somewhat more primitive compactness arguments.

\paragraph{Canonical bundles over \Teichmuller\ space} There is a smooth fiber bundle $\S\to\T$
whose fiber $\S_\sigma$ over $\sigma\in\T$ is a hyperbolic surface representing the point
$\sigma\in\T$. To make this precise, as a smooth fiber bundle we identify $\S$ with $S \cross \T$,
and we impose smoothly varying hyperbolic structures on the fibers $\S_\sigma = S \cross \sigma$,
$\sigma\in\T$, such that under the canonical homeomorphism $\S_\sigma\to S$ the hyperbolic
structure on $\S_\sigma$ represents the point $\sigma\in\T$. The action of $\MCG$ on $\T$ lifts
to a fiberwise isometric action of $\MCG$ on $\S$. Each fiber $\S_\sigma$ is therefore a
\emph{marked} hyperbolic surface, meaning that it comes equipped with an isotopy class of
homeomorphisms to $S$. The bundle $\S\to\T$ is called the \emph{canonical marked hyperbolic
surface bundle} over $\T$. 

Structures living on $\S_\sigma$ such as measured foliations or measured geodesic laminations can
be regarded as living on $S$, via the identification $\S_\sigma \approx S$. We may therefore
represent an element of $\MF$, for example, as a measured foliation on a fiber $\S_\sigma$. The
same remark will apply below, without comment, when we discuss pullback bundles of $\S$.

The \emph{canonical hyperbolic plane bundle} $\H\to\T$ is defined as the composition
$\H\to\S\to\T$ where $\H\to\S$ is the universal covering map. Each fiber $\H_\sigma$,
$\sigma\in\T$, is isometric to the hyperbolic plane, with hyperbolic structures varying smoothly
in~$\sigma$. The group $\pi_1 S$ acts as deck transformations of the covering map $\H\to\S$ is
$\pi_1 S$, acting on each fiber $\H_\sigma$ by isometric deck transformations with quotient
$\S_\sigma$. The action of $\pi_1 S$ on $\H$ extends to a fiberwise isometric action of
$\MCG(S,p)$ on $\H$, such that the covering map $\H\to\S$ is equivariant with respect to the group
homomorphism $\MCG(S,p) \to \MCG(S)$. By a result of Bers \cite{Bers:FiberSpaces}, $\H$ can be
identified with the \Teichmuller\ space of the once-punctured surface $(S,p)$, and the action of
$\MCG(S,p)$ on $\H$ is identified with the natural action of the mapping class group on
\Teichmuller\ space.

Let $T\S$ denote the tangent bundle of $\S$. Let $T_v\S$ denote the vertical sub-bundle of $T\S$,
that is, the kernel of the derivative of the fiber bundle projection $\S\to\T$, consisting of the
tangent planes to the fibers of $\S\to\T$. There exists a smoothly varying $\MCG$--equivariant
\emph{connection} on the bundle $\S\to\T$, which means a smooth sub-bundle $T_h\S$ of $T\S$ which
is complementary to $T_v\S$, that is, $T\S = T_h\S \oplus T_v\S$; see
\cite{FarbMosher:quasiconvex} for the construction of $T_h\S$. Lifting to $\H$ we have a
connection $T_h\H$ on the bundle $\H\to\T$, equivariant with respect to $\MCG(S,p)$.

\paragraph{Hyperbolic surface bundles over lines}
A \emph{closed interval} is a closed, connected subset of $\reals$, either a closed segment, a
closed ray, or the whole line. The domains of all of our paths will be closed intervals.

Given a closed interval $I \subset \reals$, a path $\gamma \from I \to \T$ is
\emph{affine} if it satisfies $d(\gamma(s),\gamma(t)) = K\abs{s-t}$ for some constant $K \ge 0$,
and $\gamma$ is \emph{piecewise affine} if there is a decomposition of $I$ into
subintervals on each of which $\gamma$ is affine. 

Given an affine path $\gamma \from I \to \T$, by pulling back the canonical marked hyperbolic
surface bundle $\S\to\T$ and its connection $T_h\S$, we obtain a marked hyperbolic surface bundle
$\S_\gamma \to I$ and a connection $T_h\S_\gamma$. This connection canonically determines a
Riemannian metric on $\S_\gamma$, as follows. Since $T_h\S_\gamma$ is \nb{1}dimensional, there is
a unique vector field $V$ on $\S_\gamma$ parallel to $T_h\S_\gamma$ such that the derivative of
the map $\S_\gamma\to I \subset \reals$ takes each vector in $V$ to the positive unit vector in
$\reals$. The fiberwise Riemannian metric on $\S_\gamma$ now extends uniquely to a Riemannian
metric on $\S_\gamma$ such that $V$ is everywhere orthogonal to the fibration and has unit
length. 

Given a piecewise affine path $\gamma\from I \to \T$, the above construction of a Riemannian
metric can be carried out over each affine subpath, and at any point $t \in I$ where two such
subpaths meet the metrics agree along the fibers, thereby producing a piecewise Riemannian metric
on $\S_\gamma$.

Given a (piecewise) affine path $\gamma \from I \to \T$, the above constructions can be carried
out on $\H_\gamma$, producing a (piecewise) Riemannian metric, equivariant with respect to $\pi_1
S$, such that the covering map $\H_\gamma\to\S_\gamma$ is local isometry.

Because our paths all have domains which are closed intervals, it follows that the induced path
metric of each Riemannian metric constructed above is a proper geodesic metric.  

A \emph{connection path} in either of the bundles $\S_\gamma\to I$, $\H_\gamma\to I$ is a
piecewise smooth section of the projection map which is everywhere tangent to the connection. By
construction, given $s,t \in I$, a path $p$ from a point in the fiber over $s$ to a
point in the fiber over $t$ has length $\le\abs{s-t}$, with equality only if $p$ is a
connection path. It follows that the min distance and the Hausdorff distance between fibers are
both equal to $\abs{s-t}$. By moving points along connection paths, for each $s,t \in I$ we have
well-defined maps $\S_s\to\S_t$, $\H_s \to \H_t$, both denoted $h_{st}$ when no confusion can
ensue. The following result gives some regularity for the maps $h_{st}$; it is closely related to
a basic fact in dynamical systems, that if $\Phi$ is a smooth flow on a compact manifold then
there is a constant $K \ge 1$ such that $\norm{D\Phi_t} \le K^{\abs{t}}$ for all $t \in \reals$.

\begin{lemma}[\cite{FarbMosher:quasiconvex}, Lemma 4.1]
\label{LemmaBilipConnection}
For each bounded set $\B\subset\Mod$ and $\rho \ge 1$ there exists $K$ such that if $\gamma\from
I \to \T$ is a $\B$--cobounded, $\rho$--lipschitz, piecewise affine path, then for each $s,t \in I$
the connection map $h_{st}$ is $K^{\abs{s-t}}$~bilipschitz.
\qed\end{lemma}

We have associated natural geometries $\S_\gamma, \H_\gamma$ to any piecewise affine path
$\gamma \from I \to \T$. When $\gamma$ is a geodesic there is another pair of natural geometries,
the singular \solv\ space $\S^\solv_\gamma$ and its universal cover $\H^\solv_\gamma$. To define
these, there is a quadratic differential $q$ such that 
$$\gamma(t) = \sigma(e^{-t} \Fol_x(q),e^{t} \Fol_y(q)), \quad t \in I.
$$
Let $\abs{dy}$ denote the transverse measure on the
horizontal foliation $f_x(q)$ and $\abs{dx}$ the transverse measure on $f_y(q)$. On the conformal
surface $\S_t \subset\S_\gamma$ we have the quadratic differential $q(e^{-t} \Fol_x(q),e^{t}
\Fol_y(q))$ with horizontal transverse measure $e^{-t}\abs{dy}$ and vertical transverse measure
$e^t \abs{dx}$. This allows us to define the \emph{singular \solv\ metric} on $\S_\gamma$ by the
formula
$$ds^2 = e^{2t} \abs{dx}^2 + e^{-2t} \abs{dy}^2 + dt^2
$$
and we denote this metric space by $\S^\solv_\gamma$. The lift of this metric to the universal
cover $\H_\gamma$ produces a metric space denoted $\H^\solv_\gamma$.

The following result says that the metric on $\H_\gamma$ is quasi-isometrically stable with
respect to perturbation of $\gamma$. Moreover, if $\gamma$ fellow travels a geodesic $\gamma'$
then the singular \solv\ geometry $\H^\solv_{\gamma'}$ serves as a model geometry.

\begin{proposition}[\cite{FarbMosher:quasiconvex}, Proposition 4.2] 
\label{PropMetricPerturbation}
For any $\rho \ge 1$, any bounded subset $\B \subset \Mod$, and any $A \ge 0$ there exists $K \ge
1$, $C \ge 0$ such that the following holds. If $\gamma,\gamma' \from I \to\T$ are two
$\rho$--lipschitz, $\B$--cobounded, piecewise affine paths defined on a closed interval $I$, and if
$d(\gamma(s),\gamma'(s)) \le A$ for all $s\in I$, then there exists a map
$\S_\gamma\to\S_{\gamma'}$ taking each fiber $\S_{\gamma(t)}$ to the fiber $\S_{\gamma'(t)}$ by a
homeomorphism in the correct isotopy class, such that any lifted map $\H_\gamma \to \H_{\gamma'}$
is a $K,C$--quasi-isometry.

If $\gamma'$ is a geodesic, the same is true with
$\S_{\gamma'}$, $\H_{\gamma'}$ replaced by the singular \solv\ spaces $\S^\solv_{\gamma'}$,
$\H^\solv_{\gamma'}$.
\end{proposition}

\paragraph{Remark} Given any $\rho$--lipschitz path $\gamma\from I\to\T$, the bundles $\S\to\T$,
$\H\to\T$ can be pulled back to define bundles $\S_\gamma \to I$, $\H_\gamma \to I$. Despite the
paucity of smoothness, one can extend the fiberwise hyperbolic metrics on these bundles to
measurable Riemannian metrics which then determine proper geodesic metrics in a canonical manner,
and hence we would associate a geometry $\H_\gamma$ to $\gamma$. However, it is easier to proceed
by approximating $\gamma$ with a $\rho$--lipschitz piecewise affine path, and to apply
Proposition~\ref{PropMetricPerturbation} to show that the resulting geometry on $\H_\gamma$ is
well-defined up to quasi-isometry. This allows us to reduce the proof of
Theorem~\ref{TheoremStableTeichGeod} to the case where $\gamma$ is piecewise affine, a technical
simplification.

\section{Proof of Theorem \ref{TheoremStableTeichGeod}}
\label{SectionTeichGeodProof}

For the proof we fix the closed, oriented surface $S$ of genus $\ge 2$, with \Teichmuller\ space
$\T$, mapping class group $\MCG$, measured foliation space $\MF$, etc.

\subsection{Setting up the proof}

Throughout the proof we fix a compact subset $\B \subset \Mod$, and numbers $\rho \ge 1$,
$\delta \ge 0$. 

Given $\gamma \from I \to \T$ a $\B$--cobounded, $\rho$--Lipschitz, path,
Proposition~\ref{PropMetricPerturbation} says that if we perturb $\gamma$ to be a piecewise
affine path, then the large scale geometry of $\H_\gamma$ is well-defined up to quasi-isometry. 
The truth of the hypothesis of Theorem~\ref{TheoremStableTeichGeod} is therefore unaffected by
perturbation, as is the conclusion, with uniform changes in all constants depending only on the
size of the perturbation. We shall choose a particular perturbation which will be technically
useful in what follows, particularly in Section~\ref{SectionCompactnessProp}.

A \emph{$\Z$--piecewise affine path} is a path $\gamma \from I \to \T$, defined on a closed
interval~$I$ whose endpoints, if any, are in $\Z$, such that $\gamma$ is affine on each
subinterval $[n,n+1]$ with $n,n+1 \in I \intersect \Z$. We shall often denote
$J=I\intersect\Z$. If $\gamma$ is a $\rho$--lipschitz path defined on an interval $I$, we can
perturb $\gamma$ to be $\Z$--piecewise affine as follows: first remove a bit from each finite end
of $I$ of length less than $1$ so that the endpoints are in $\Z$; then replace $\gamma
\restrict [n,n+1]$ by an affine path with the same endpoints whenever $n,n+1 \in J$; note
that $d(\gamma(t),\gamma'(t)) \le \rho$ for all $t \in I$. 

Given a $\B$--cobounded, $\rho$--lipschitz, $\Z$--piecewise affine path $\gamma\from I \to \T$ whose
canonical hyperbolic plane bundle $\H_\gamma \to I$ is $\delta$--hyperbolic, our goal is to
construct a \Teichmuller\ geodesic $g$, sharing any endpoints with $\gamma$, and to show that
$\gamma$ and $g$ are at bounded Hausdorff distance and $\gamma$ is a quasigeodesic. 

\paragraph{Motivation} The construction of $g$ is inspired by ending laminations methods
(\cite{Thurston:GeomTop}, chapter 9; \cite{Bonahon:ends}), uniform foliations methods
\cite{Thurston:FoliationsAndCircles}, and the flaring concepts from
\cite{BestvinaFeighn:combination}. 

In the case of a line $\gamma\from\reals\to\T$, the idea for finding $g$ is to keep in mind an
analogy between the hyperbolic surface bundle $\S_\gamma\to\reals$ and a doubly degenerate
hyperbolic structure on $S \cross \reals$. The hyperbolic surfaces $\S_t\subset\S_\gamma$ ($t \in
\reals$) approach both ends of $\S_\gamma$ as $t \to\pm\infinity$ and so, following the analogy,
$\S_\gamma$ is ``geometrically tame'' in the sense of \cite{Bonahon:ends} and
\cite{Thurston:GeomTop} chapter~9. This suggests that we produce an ending lamination in $\PMF$
for each of the two ends. This pair of laminations will jointly fill the surface and so will
determine a geodesic $g$ in $\T$. 

Despite the analogy, our construction of ending laminations for $\S_\gamma$ is entirely
self-contained and new. The construction is inspired by uniform foliations methods
\cite{Thurston:FoliationsAndCircles}, in which large-scale geometry of a foliation is used to
determine laminations on leaves of that foliation. The main new idea is that, in the presence of
$\delta$--hyperbolicity, laminations can be constructed using flaring concepts from
\cite{BestvinaFeighn:combination}. The only prerequisites for the construction are basic facts
about measured geodesic laminations and about flaring.

\subsection{Flaring}
\label{SectionFlaring}

In order to get the proof off the ground, the key observation needed is that hyperbolicity of
$\H_\gamma$ is equivalent to the ``rectangles flare'' property introduced by Bestvina and Feighn
\cite{BestvinaFeighn:combination}. Sufficiency of the rectangles flare property was proved
by Bestvina and Feighn, and necessity was proved by Gersten \cite{Gersten:cohomological}. In
\cite{FarbMosher:quasiconvex} the rectangles flare property is recast in a manner which is
followed here.

Consider a sequence of positive numbers $r_j$ $(j \in J)$ indexed by a subinterval $J$ of the
integers $\Z$, which can be finite, half-infinite, or all of $\Z$. Given $\kappa > 1$,
$n \in \Z_+$, $A \ge 0$, the \emph{$\kappa,n,A$--flaring property} says that if $j-n,j,j+n \in J$,
and if $r_j > A$, then $\max\{r_{j-n},r_{j+n}\} \ge \kappa r_j$. Given $L \ge 1$, the
\emph{$L$--lipschitz growth condition} says that for every $j,k \in J$ we have $r_j/r_k
\le L^{\abs{j-k}}$ for $j,k \in J$; equivalently, if $\abs{j-k}=1$ then $r_j \le L r_k$. Given $L
\ge 1$, $D \ge 0$, the \emph{$(L,D)$ coarse lipschitz growth condition} says that if
$\abs{j-k}=1$ then $r_j \le L r_k + D$.

The flaring property and the coarse lipschitz growth condition work together. The
$\kappa,n,A$ flaring property alone only controls behavior on arithmetic subsequences of the form
$j_0 + kn$, e.g.\ if $r_{j_0} > A$ then there exists $\epsilon=\pm 1$ 
such that 
\begin{equation}
r_{j_0+\epsilon k n} > r_j\kappa^k \quad\text{for all}\quad k\ge 1.
\label{EquationGrowth}
\end{equation}
But in company with, say, the $L$--lipschitz growth
condition, growth of disjoint arithmetic subsequences is conjoined, in that 
\begin{equation}
r_{j_0 + \epsilon m} > r_j L^{-n}\kappa'{}^m, \quad\text{for all}\quad m \ge 1.
\label{EquationFullGrowth}
\end{equation}
where $\kappa'=\kappa^{1/n}$. A similar result holds for a coarse lipschitz growth condition.

The number $A$ is called the \emph{flaring threshold}. A $\kappa,n,A$ flaring sequence can stay
$\le A$ on a subinterval of arbitrary length, but as we have seen, at any point where the
sequence gets above the flaring threshold then exponential growth is inexorable in at least one of
the two directions. In particular, if the flaring threshold is zero then exponential growth is
everywhere.

Consider now a cobounded, lipschitz, piecewise affine path $\gamma \from I \to \T$ and the
hyperbolic plane bundle $\H_\gamma\to I$. Given $s \in I$ let $d_s$ denote the distance function
on the fiber $\H_s$. A \emph{$\lambda$--quasihorizontal path} in $\H_\gamma$, $\lambda\ge 1$, is
a section $\alpha \from I \to \H_\gamma$ of the bundle projection which is $\lambda$--lipschitz. A
\nb{1}quasihorizontal path is the same thing as a connection path. For any pair of
$\lambda$--quasihorizontal paths $\alpha,\beta \from I \to \H_\gamma$, the sequence of distances 
$$d_j(\alpha(j),\beta(j)), \quad j \in J=I\intersect\Z
$$
automatically satisfies an $(L,D)$ coarse Lipschitz growth condition, where $L=K$ is the constant
given by Lemma~\ref{LemmaBilipConnection} and hence depends only on the coboundedness and the
lipschitz constant of $\gamma$, and where $D=2K(\lambda+1)$. To see this, given $j,k \in J$ with
$\abs{j-k}=1$, let $a=\alpha(j)$, $a'=\alpha(k)$, $a''=h_{jk}(a)$, $b=\beta(j)$, $b'=\beta(k)$,
$b''=h_{jk}(b)$. The points $a',a''$ are connected by a path of length $\le\lambda+1$
staying between $\H_j$ and $\H_k$, consisting of a segment of $\alpha$ from $a'$ to $a$ and a
connection path from $a$ to $a''$; projection of this path onto $\H_k$ is $K$--lipschitz, and
similarly for the $b$'s, and so:
\begin{align*}
d_k(a',a'') &\le K(\lambda+1) \\
d_k(b',b'') &\le K(\lambda+1) \\
d_k(a',b') &\le d_k(a',a'') + d_k(a'',b'') + d_k(b'',b') \\
   &\le K(\lambda+1) + K d_j(a,b) + K(\lambda+1) 
\end{align*}
Given constants $\kappa \ge 1$, $n \in \Z_+$, and a function $A(\lambda) \ge 0$ defined for
$\lambda \ge 1$, we say that $\H_\gamma$ satisfies the \emph{$\kappa,n,A(\lambda)$ horizontal
flaring property} if for any $\lambda \ge 1$ and any $\lambda$--quasihorizontal paths
$\alpha,\beta \from I \to \H_\gamma$, the sequence of distances
$$d_j(\alpha(j),\beta(j)), \quad j \in J=I\intersect\Z
$$
satisfies $\kappa,n,A(\lambda)$ flaring. 

The following result is an almost immediate consequence of Gersten's theorem
\cite{Gersten:cohomological} which gives the converse to the Bestvina--Feighn combination
theorem; the interface with Gersten's theorem is explained in \cite{FarbMosher:quasiconvex}.
See also Lemma~5.2 of \cite{FarbMosher:quasiconvex} for an alternative proof following the lines
of the well-known fact that in a hyperbolic metric space, geodesic rays satisfy an exponential
divergence property \cite{Cannon:TheoryHyp}.

\begin{proposition} Given a $\rho$--lipschitz, $\B$--cobounded map $\gamma \from I \to \T$, if
$\H_\gamma$ is $\delta$--hyperbolic, then $\H_\gamma$ satisfies a horizontal
flaring property, with flaring data $\kappa,n,A(\lambda)$ dependent only $\rho,\B,\delta$.
\qed
\end{proposition}

In the context of the proof of Theorem~\ref{TheoremStableTeichGeod}, for each $\mu \in \MF$, each
marked hyperbolic surface $\S_t$ in the bundle $\S_\gamma$ has a measured geodesic lamination
$\mu_t$ in the class $\mu$. Denote $\Length_t(\mu) = \Length_{\gamma(t)}(\mu_t)$. The
key to our proof of Theorem~\ref{TheoremStableTeichGeod} is to study how the function
$\Length_t(\mu)$ varies in $t$ and in $\mu$. What gets us off the ground is
Lemma~\ref{LemmaLengthFlares} which shows that for each $\mu\in \MF$ the function
$t\mapsto\Length_t(\mu)$ flares, with uniform flaring data independent of
$\mu$, and with a flaring threshold of zero; intuitively this follows from the flaring of
$\H_\gamma$. The technical step of verifying uniform flaring data depends on the
construction of a flow preserving connection on the leafwise geodesic flow bundle of
$\S_\gamma$, to which we now turn.

\subsection{Connections on geodesic flow bundles} 

Let $\GFL(F)$ denote the geodesic flow of a hyperbolic surface $F$, defined on the unit
tangent bundle $T^1 F$ of $F$. The foliation of $\GFL(F)$ by flow lines, being invariant under
the antipodal map $\vec v \mapsto -\vec v$ of $T^1 F$, descends to a foliation on the tangent line
bundle $\P T^1 F = T^1 F / \pm$, the \emph{geodesic foliation} $\GF(F)$. The Liouville metric on
$T^1 F$ descends to a metric on $\P T^1 F$ also called the Liouville metric.

Starting from the canonical marked hyperbolic surface bundle $\S\to\T$, by taking the geodesic
flow on each fiber we obtain the \emph{fiberwise geodesic flow bundle} $\GFL(\S)\to\T$ whose fiber
over $\sigma\in\T$ is $\GFL(\S_\sigma)$, and we similarly obtain the fiberwise geodesic foliation
bundle $\GF(\S) \to \T$. Note that the geodesic flows on the fibers of $\GFL(\S)$ fit together to
form a smooth flow on $\GFL(\S)$, and similarly for the geodesic foliations on fibers of
$\GF(\S)$; smoothness follows from the fact that the leafwise hyperbolic metrics on $\S \to \T$
vary smoothly, together with the fact that the geodesic differential equation has coefficients
and also solutions varying smoothly with the metric. We have a fiber bundle $\GFL(\S) \to \S$
whose fiber over $x \in \S_\sigma \subset \S$ is the unit tangent space
$T^1_x\S_\sigma$, and similarly a fiber bundle $\GF(\S) \to \S$ whose fiber is the space
of tangent lines $\P T^1_x \S_\sigma$.

For any piecewise affine path $\gamma \from I \to \T$, we have pullback bundles
$\GFL(\S_\gamma) \to I$ and $\GF(\S_\gamma) \to I$, whose fibers over $t \in I$ are $\GFL(\S_{t})$
and $\GF(\S_{t})$ respectively. 

In the case where $\gamma$ is affine, by a \emph{connection} on the bundle $\GF(\S_\gamma) \to I$
we mean a connection which preserves the geodesic foliations, that is: a \nb{1}dimensional
sub-bundle of the tangent bundle of $\GF(\S_\gamma)$ which is complementary to the vertical
sub-bundle of $\GF(\S_\gamma) \to I$, such that for any $s,t \in I$, the map
$H_{st} \from \GF(\S_s) \to \GF(\S_t)$ obtained by moving points along connection paths takes
leaves of $\GF(\S_s)$ to leaves of $\GF(\S_t)$. The \emph{connection flow} $\Phi$ on
$\GF(\S_\gamma)$ is defined by $\Phi_r(l) = H_{s,s+r}(l)$ whenever $l \in \GF(\S_s)$. 

In the general case where $\gamma$ is only piecewise affine, the connection flows defined over
the intervals where $\gamma$ is affine piece together to define a connection flow $\Phi$ on all
of $\GF(\S_\gamma)$, with corresponding connection maps $H_{st}(l) = \Phi_{t-s}(l)$ for all $s,t
\in I$, $l \in \GF(\S_s)$. 

A connection on $\GF(\S_\gamma)$ is \emph{leafwise $L$--bilipschitz} if $H_{st}$ restricts to a
$L^{\abs{s-t}}$--bilipschitz homeomorphism between leaves of the geodesic foliation.

All the above concepts apply as well to the fiberwise geodesic flow bundle and geodesic foliation
bundle $\GFL(\H)\to\T$, $\GF(\H)\to\T$, and any associated pullback bundles.

\begin{lemma} 
\label{LemmaGFConnection}
For each bounded set $\B\subset\Mod$ and each $\rho \ge 1$ there exists $L \ge 1$ and $\lambda
\ge 1$ such that if $\gamma\from I \to \T$ is a $\B$--cobounded, $\rho$--lipschitz, piecewise
affine path, then there is a leafwise $L$--bilipschitz connection on the bundle $\GF(\S_\gamma) \to
I$, and each connection line in $\GF(\S_\gamma)$ projects to a path in $\S_\gamma$ which is
$\lambda$--quasihorizontal. By lifting that we obtain a $\pi_1(S)$--equivariant leafwise
$L$--bilipschitz connection on $\GF(\H_\gamma)$ whose connection lines project to
$\lambda$--quasihorizontal paths in $\H_\gamma$. Moreover the connection satisfies the following
uniform continuity condition: for every bounded $\B\subset\Mod$, $\rho \ge 1$, $M > 0$, and
$\epsilon>0$, there exists $\delta>0$ such that if $\gamma \from I \to \T$ is
$\B$--cobounded, $\rho$--lipschitz, and piecewise affine, if $s,t \in I$ with $\abs{s-t}  < M$, and
if $l,m \in \GF(\S_s)$ with $d(l,m) < \delta$, then $d(H_{st}(l),H_{st}(m)) < \epsilon$. 
\end{lemma}

I am grateful to Amie Wilkinson for the proof of this lemma, particularly for explaining how to
apply partially hyperbolic dynamics.

\begin{proof} The intuition behind the proof is that the geodesic flows on $\S_s$ and $\S_t$ are
topologically conjugate; to put it another way, for each geodesic $\ell$ on $\S_s$, the
connection on $\S_\gamma$ moves $\ell$ to a $K^{\abs{s-t}}$--bilipschitz path in $\S_t$
and hence that path is close to a geodesic in $\S_t$. In order to carry this out uniformly up
in the geodesic foliation bundle $\GF(\S_\gamma)$ we shall apply structural stability tools
from the theory of hyperbolic dynamical systems, as encapsulated in the
Sublemma~\ref{SublemmaGeodesicFoliation}.

A $k$--dimensional foliation of a subset of $\R^n$ has \emph{uniformly smooth leaves} if the
leaves are defined locally by immersions from open subsets of $\R^k$ into $\R^n$ such that for
each $r \ge 0$ the partial derivatives up to order $r$ are uniformly bounded away from zero and
from infinity. A foliation of a smooth manifold $M$ has \emph{locally uniformly smooth leaves} if
$M$ is covered by coordinate charts in each of which the leaves are uniformly smooth. The
property of \emph{(locally) uniformly $C^r$ leaves} is similarly defined by omitting the words
``for each $r \ge 0$''.

Let $D$ be the dimension of $\T$.

\begin{sublemma} 
\label{SublemmaGeodesicFoliation}
There exists a unique, $D+1$ dimensional foliation $\G$ of $\GFL(\S)$ with locally uniformly
smooth leaves such that $\G$ is transverse to the fibers of $\GFL(\S) \to \T$, and the foliation
of $\GFL(\S)$ obtained by intersecting $\G$ with the fibers of $\GFL(\S)\to\T$ is identical to the
foliation by geodesic flow lines.
\end{sublemma}

Before proving the claim we apply it to prove Lemma~\ref{LemmaGFConnection}. 

Since the conclusion of the lemma is local, it suffices to prove it when $\gamma$ is a
$\rho$--lipschitz, $\B$--cobounded affine arc $\gamma \from [0,1]\to\T$.

Choose a compact set $\A \subset \T$ such that any $\rho$--lipschitz, $\B$--cobounded affine arc
$[0,1] \to \T$ may be translated by $\Isom(\T)$ to lie in $\A$. Let $C_\rho(\A)$ be the space
of all $\rho$--lipschitz affine arcs $[0,1] \to \A$. The conclusion of the lemma is invariant under
the action of $\MCG$ and so we may assume $\gamma \in C_\rho(\A)$. By enlarging $\A$ we may assume
that $\A$ is a smooth codimension~0 submanifold of~$\T$. By restricting $\S$ to $\A$ we obtain a
hyperbolic surface bundle $\S_\A$ and its geodesic flow bundle $\GFL(\S_\A)$; fix a smooth
Riemannian metric on $\GFL(\S_\A)$. Let $\G_\A$ be the restriction of the foliation $\G$ to
$\GFL(\S_\A)$. For each $\gamma \in C_\rho(\A)$, the foliation $\G_\A$ restricts to a
2--dimensional foliation $\G_\gamma$ of $\GFL(\S_\gamma)$ with uniformly smooth leaves, transverse
to the fibers of $\GFL(\S_\gamma) \to [0,1]$. Also, the Riemannian metric on $\GFL(\S_\A)$
restricts to a smooth Riemannian metric on $\GFL(\S_\gamma)$. There is a unique vector field
$V_\gamma$ on $\GFL(\S_\gamma)$ which is tangent to $\G_\gamma$ and is perpendicular to the
geodesic flow lines, such that each $v \in V_\gamma$ projects to a positive unit tangent vector
in $[0,1] = \domain(\gamma)$. By uniqueness the foliation $\G$ is invariant under the antipodal
map on $\GFL(\S)$, implying that $\G_\gamma$ is antipode invariant on
$\GFL(\S_\gamma)$. Assuming as we may that the Riemannian metric on $\GFL(\S_\A)$ is also
antipode-invariant, it follows that $V_\gamma$ is antipode-invariant, and so descends to a
vector field on $\GF(\S_\gamma)$. This vector field spans the desired connection on
$\GF(\S_\gamma)$. 

Note that the connection on $\GF(\S_\gamma)$ is uniformly smooth along leaves of $\G_\gamma$, and
as $\gamma$ varies over $C_\rho(\A)$ the connection varies continuously; this follows from
leafwise uniform smoothness of $\G$.

Let $\psi_{\gamma,t}$ be the connection flow on $\GF(\S_\gamma)$, which has connection maps
$H_{st} \colon$\break$\GF(\S_s) \to \GF(\S_t)$, that is, $H_{st}(\ell) = \psi_{\gamma,t-s}(\ell)$ for
$\ell \in \GF(\S_s)$. To prove that $H_{st}$ is $L^{\abs{s-t}}$ bilipschitz, it
suffices by a standard result in differential equations to prove that
$$S(\gamma) = \abs{ \frac{d}{dt}  \norm{ D\psi_{\gamma,t} } }_{t=0} \le \log(L).
$$
Note that $S(\gamma)$ is continuous as a function of $\gamma \in C_\rho(\A)$. Since
$C_\rho(\A)$ is compact, $S(\gamma)$ has a finite upper bound $l \ge 0$, and so $L=e^l$ is the
desired bilipschitz constant.

Further compactness arguments show that the connection lines project to
$\lambda$--quasihorizontal lines in $\S_\gamma$ for uniform $\lambda$, and that the uniform
continuity clause holds.
\end{proof}

\begin{proof}[Proof of Sublemma \ref{SublemmaGeodesicFoliation}]
Uniqueness of $\G$ follows because, for any $\sigma \in \T$ and any closed geodesic $c$ in
$\S_\sigma$, and for any $\tau\in\T$, the leaf of $\G$ containing $c$ must contain the closed
geodesic in $\S_\tau$ that is in the isotopy class of $c$ with respect to the canonical
homeomorphism $\S_\sigma \approx \S_\tau$. The non-simply connected leaves of $\G$ are therefore
determined, but they are dense in $\GFL(\S)$ and so $\G$ is determined.

Existence of $\G$ is a purely local phenomenon, because if we have open subsets $U,V \subset \T$
and foliations $\G_U, \G_V$ on $\GFL(\S_U), \GFL(\S_V)$ respectively, so that $\G_U,\G_V$ each
satisfy the conclusions of the sublemma, then the uniqueness argument above can be applied locally
to show that the restrictions to $\GFL(\S_{U \intersect V})$ of $\G_U,\G_V$ are identical, and so
$\G_U,\G_V$ can be pasted together over $U \union V$ to give a foliation of $\GFL(\S_{U \union
V})$ satisfying the conclusions of the sublemma. Arguing similarly with respect to some locally
finite open cover of $\T$ allows one to construct $\G$.

For each $\sigma \in \T$ it therefore suffices to find some closed ball $B$ in $\T$ around
$\sigma$ and a $D+1$ dimensional foliation $\G_B$ of $\GFL(\S_B)$ satisfying the conclusions of
the sublemma, namely, that $\G_B$ has uniformly smooth leaves and is transverse to the fibers of
$\GFL(\S_B) \to B$, and the intersection of $\G_B$ with each fiber of is the foliation of that
fiber by geodesic flow lines.

We review some elements of partially hyperbolic dynamical systems from \cite{HirschPughShub}.
Let $M$ be a smooth compact Riemannian manifold. A $C^r$--flow $\phi \from M \cross \reals \to
M$ is \emph{$r$--normally hyperbolic} at a foliation $\F$ if there is a splitting $TM = E^u
\oplus T\F \oplus E^s$, invariant under the flow $\phi$, and there exists $a>0$, such that
for each $t>0$ we have
\begin{itemize}
\item $\norm{T\phi_t \restrict E^s} < e^{-at}$,
\item $\norm{T\phi_{-t}\restrict E^u} < e^{-at}$ , 
\item $\norm{T\phi_{-t} \restrict T\F}^{r} \cdot \norm{T\phi_t \restrict E^s} < e^{-at}$,
\item $\norm{T\phi_t \restrict T\F}^r \cdot \norm{T\phi_{-t} \restrict E^u} < e^{-at}$.
\end{itemize}

We need the following results from \cite{HirschPughShub}:

\begin{theorem}
\label{TheoremStructuralStability}
Suppose that the flow $\phi$ is $r$--normally hyperbolic at $\F$.  
If the foliation $\F$ is $C^1$, then for every $C^r$ flow $\psi$, sufficiently $C^1$--close to
$\phi$, there exists a foliation $\G$ such that $\psi$ is $r$--normally hyperbolic at $\G$, and
such that the dimensions of corresponding summands in the splittings of $TM$ associated to $\F$
and to $\G$ are identical.
\end{theorem}

\begin{theorem}
\label{TheoremSmoothing}
If the flow $\psi$ is $r$--normally hyperbolic at the foliation $\G$, then the
leaves of $\G$ are uniformly $C^r$.
\end{theorem}

Given a smooth closed ball $B$ in $\T$, consider the $D+3$ dimensional manifold $M=T^1 \S_B$, the
fiberwise unit tangent bundle of $\S_B$, whose fiber over the point $x\in\S_\tau\subset\S_B$,
where $\tau\in B$, is $T^1_x(\S_\tau)$. There is also a fibration of $T^1\S_B$ over $B$, whose
fiber over $\tau \in B$ is $T^1 \S_\tau$. As a manifold, the space $T^1\S_B$ is identified with
the underlying space of the geodesic flow bundle $\GFL(\S_B)$. Let $\psi$ be the fiberwise
geodesic flow on $T^1 \S_B$ (i.e.\ the flow on $\GFL(\S_B)$). 

We claim that if $\G$ is a \nb{codimension}2 foliation of $T^1 \S_B$ at which $\psi$ is
$r$--normally hyperbolic, so that the bundles $E^s$, $E^u$ are \nb{1}dimensional, then
$\G$ is transverse to the fibers of the fibration $T^1 \S_B \to B$, and the intersection of $\G$
with each fiber $T^1\S_\tau$, $\tau\in B$, is the geodesic flow of the hyperbolic surface
$\S_\tau$. Uniform smoothness of leaves of $\G$ follows from Theorem~\ref{TheoremSmoothing},
thereby proving Sublemma~\ref{SublemmaGeodesicFoliation}. 

To prove the claim, note that the \nb{1}dimensional line bundle $T\psi$ is a sub-bundle of
$T\G$, and moreover $T\psi$ is tangent to each the fibers $T^1\S_\tau$. Each of these fibers is
\nb{3}dimensional, and so transversality of $\G$ to these fibers will follow by proving that in
the splitting $TM = E^u \oplus T\G \oplus E^s$, the sub-bundle $E^u \oplus T\psi \oplus E^s$ is
identical to the fiberwise tangent bundle of the fibration $T^1 \S_B \to B$. This will follow in
turn by proving that each of the \nb{1}dimensional bundles $E^u$, $E^s$ is tangent to the fibers
$T^1\S_\tau$. The restriction of $\psi$ to each fiber $T^1 \S_\tau$ is the geodesic flow, which is
known to be an Anosov flow with \nb{1}dimensional stable and unstable bundles. Moreover, for any
vector $v \in TM$ which is not tangent to a fiber $T^1 \S_\tau$, the component of $v$
transverse to the fibers is preserved in norm by the flow $\psi$, with respect to a Riemannian
metric on $M$ that assigns constant distance from any point in one fiber to any other fiber, and
so $v$ cannot be in the stable or unstable bundle of the splitting $E^s \oplus T\G \oplus E^u$. 
It follows that $E^u$ is the same as the Anosov unstable bundle of the fiberwise geodesic flow on
$T^1 \S_B$, and $E^s$ is the same as the Anosov stable bundle, and so $E^s$, $E^u$ are indeed
tangent to the fibers $T^1 \S_\tau$. And since $T\psi$ is a sub-bundle of $T\G$ and $E^s$, $E^u$
are independent of $T\G$, it follows that the restriction of $\G$ to each fiber $T^1 \S_\tau$ is
indeed the geodesic flow.

We have therefore reduced the proof of Sublemma~\ref{SublemmaGeodesicFoliation} to the
construction of a foliation $\G$ on $T^1 \S_B$ at which the fiberwise geodesic flow is
$r$--normally hyperbolic, for some open neighborhood $B$ of any point $\sigma \in \T$. We carry
out this construction using Theorem~\ref{TheoremStructuralStability}.

Fix a point $\sigma$ of $\T$. To start with let $B$ be any smooth closed ball in $\T$ whose
interior contains $\sigma$. Take $M=T^1 \S_B$ as above. Pick a diffeomorphism
$\S_{B}\homeo\S_\sigma \cross B$ respecting projection to $B$. This induces a diffeomorphism
$\Theta \from T^1 \S_B \to T^1\S_\sigma\cross B$.  The geodesic flow on $T^1 \S_\sigma \cross B$
pulls back via $\Theta$ to a flow $\phi$ on $T^1 \S_B$. Also, there is a $D+1$ dimensional
foliation of $T^1\S_\sigma \cross B$ which is the product of the geodesic flow on $T^1 \S_\sigma$
with $B$; pulling this foliation back via $\Theta$ we obtain a foliation of $T^1 \S_B$ denoted
$\F$. Noting that $\norm{T\phi_t\restrict T\F} = 1$ for all $t$, from the fact that the geodesic
flow on $\S_\sigma$ is Anosov it follows that $\phi$ is $r$--normally hyperbolic at $\F$
for all $r$, with \nb{1}dimensional stable and unstable bundles $E^s$, $E^u$.

For the flow $\psi$ on $T^1 \S_B$ we would like to take the leafwise geodesic flow. However,
we have no control on the $C^1$ distance between $\phi$ and $\psi$ as required to apply
Theorem~\ref{TheoremStructuralStability}. To fix this, we want to ``choose the ball $B$ to be
sufficiently small'', but we must do this in a way that does not change the domain manifold
$M = T^1 \S_B$. Choose a diffeomorphism between $(B,\sigma)$ and the unit ball in Euclidean
space centered at the origin. With respect to this diffeomorphism let $(0,1) \cross B  \to B$
denote scalar multiplication, and so $s \cdot B$ corresponds to the ball of radius $s$ in
Euclidean space. As $s$ approaches zero, lift the maps $B \to s \cdot B$ to a smooth family
of diffeomorphisms $T^1 \S_B \to T^1 \S_{s \cdot B}$ so that for each $b \in B$ the
diffeomorphisms $T^1 \S_b \to T^1 \S_{s\cdot b}$ converge uniformly to $\Theta \from T^1\S_b
\to T^1\S_\sigma$ as $s \to 0$. Pulling back the leafwise geodesic flow on $T^1\S_{s \cdot B}$
we obtain a smooth family of flows $\psi_s$ on $T^1\S_B$ converging to the flow $\phi$ in the
$C^1$ topology, as $s \to 0$. 

We may now apply Theorem~\ref{TheoremStructuralStability} to $\psi_s$ for $s$ sufficiently
close to $0$, to obtain a foliation $\G_s$ on $T^1\S_B$ at which $\psi_s$ is $r$--normally
hyperbolic. Pulling back to the ball $B' = s \cdot B$ around $\sigma$ we have constructed a
foliation $\G_{B'}$ at which the leafwise geodesic flow is $r$--normally hyperbolic.
\end{proof}

\subsection{Flaring of geodesic laminations}
\label{SectionFlaringLaminations}

Until further notice we shall fix a $\Z$--piecewise affine path $\gamma \from I \to \T$ which is
$\B$--cobounded and $\rho$--lipschitz, such that $\H_\gamma$ is $\delta$--hyperbolic. Let $J = I
\intersect \Z$. 

For each $\mu \in \MF$, let $\mu_t$ denote the measured geodesic lamination on $\S_t$
representing $\mu$. For $i \in J$ let $\Length_i(\mu)$ be the length of $\mu_i$ in the hyperbolic
surface $\S_i$. The following lemma says that the sequence $\Length_i(\mu)$ flares with uniform
flaring constants, and with a flaring threshold of \emph{zero}. 

\begin{lemma}[Length flares] 
\label{LemmaLengthFlares}
There exist constants $L \ge 1$, $\kappa > 1$, $n \in \Z_+$ depending only on $\B,\rho,\delta$
such that the following holds. For any $\mu\in\MF$, the sequence $i \mapsto \Length_i(\mu)$, $(i
\in J)$, satisfies the $L$--lipschitz, $(\kappa,n,0)$ flaring property.
\end{lemma}

\begin{proof} From Lemma~\ref{LemmaBilipConnection} we have a $K$--bilipschitz connection on
$\S_\gamma$, and from Lemma~\ref{LemmaGFConnection} we have an $L$--bilipschitz connection on
$\GF(\S_\gamma)$ whose connection lines project to $\lambda$--quasihorizontal lines in $\S_\gamma$.
By lifting we obtain similar connections in $\H_\gamma$ and $\GF(\H_\gamma)$. The constants
$K,L,\lambda$ depend only on $\B$,~$\rho$.

Define the \emph{suspension} of $\mu$ to be the following
\nb{2}dimensional measured lamination in $\S_\gamma$:
$$\Susp(\mu) = \Union_{t \in \reals} \mu_t
$$
By Theorems~\ref{TheoremStructuralStability} and \ref{TheoremSmoothing}, it follows that each
leaf of $\Susp(\mu)$ is piecewise smooth, being uniformly smooth over each affine segment of
$\gamma$. Restricting the projection $\S_\gamma \to \reals$ to $\Susp(\mu)$, we may think of
$\Susp(\mu)$ as a $\mu$--bundle over~$I$.

By restricting the connection on $\GF(\S_\gamma)$, we obtain a connection on
$\Susp(\mu)$ whose connection lines respect the leaves of $\Susp(\mu)$ and are transverse to the
fibers $\mu_t$. By lifting to the universal cover $\H_\gamma$ of $\S_\gamma$ we obtain the
suspension $\Susp(\wt\mu)$ of $\wt\mu$, whose fiber in $\H_t$ is $\wt\mu_t$, and we obtain a
connection on $\Susp(\wt\mu)$. We will use $h_{st}$ to denote either of the connection maps $\mu_s
\to \mu_t$ or $\wt\mu_s \to \wt\mu_t$; the context should make the meaning clear. Note in
particular that $h_{st}$ is $L^{\abs{s-t}}$ bilipschitz from leaves of $\wt\mu_s$ to leaves of
$\wt\mu_t$, preserves transverse measure, and has connection lines which are $\lambda$
quasihorizontal in~$\H_\gamma$. 

Since $\H_\gamma$ is $\delta$--hyperbolic, it satisfies horizontal flaring with data depending only
on $\B$, $\rho$, $\delta$, and so we have:

\begin{lemma}[Uniform Flaring]
\label{LemmaUniformFlaring} 
There exist constants $(\kappa,n,A)$ depending only on $\B$, $\rho$, $\delta$ with the following
property. For any $\mu \in \MF$, and for any two connection lines $\alpha,\alpha'$ of
$\Susp(\wt\mu)$, the sequence of distances $d_{\H_i}(\alpha(i),\alpha'(i))$ satisfy the
$L$--lipschitz, $(\kappa,n,A)$ flaring property. In particular, for any $\mu \in \MF$, any $s \in
\reals$, and any leaf segment $\ell$ of $\mu_s$, the sequence of lengths
$\Length_{s+i}(h_{s,s+i}(\ell))$ satisfies the $L$--lipschitz, $(\kappa,n,A)$ flaring property.
\qed\end{lemma}

\paragraph{Remark} Connection lines of $\Susp(\wt\mu)$ are only quasihorizontal, not horizontal,
so they do \emph{not} necessarily coincide with connection lines in $\H_\gamma$.
\medskip

Continuing the proof of Lemma~\ref{LemmaLengthFlares}, recall that:
$$\Length_t(\mu) = \int_{\mu_t} d\mu_t = \int_{\mu_t} d\mu_t^\perp \cross d\mu_t^\parallel
$$
Using this formula we can express the length change from $\Length_s(\mu)$ to $\Length_t(\mu)$
as an integral of a derivative. To be precise, as a measurable function on $\mu_s$ we
have a well-defined Radon--Nykodym derivative:
$$h'_{st} = \frac{h_{st*}^\inv (d\mu_t)}{d\mu_s} = \frac{h_{st*}^\inv
(d\mu^\parallel_t)}{d\mu^\parallel_s}
$$
It follows that:
$$\Length_t(\mu) = \int_{\mu_t} d\mu_t = \int_{\mu_s} h_{st}' \, d\mu_s
$$
For each $x \in \mu_s$ let $I_a(x)$ be the lamination segment of
length $a$ centered on~$x$. By applying Fubini's theorem and using a change of
variables, we get:
\begin{align*}
\Length_t(\mu) &= \int_{\mu_s} \left( \int_{\tau \in I_a(x)} \frac{h_{st}'(\tau)}{a}
\, d\tau
\right) d\mu_s(x) \\
\intertext{where $d\tau$ is simply leafwise Lebesgue measure, that is, $d\tau =
d\mu_s^\parallel$. We rewrite this as:} 
\Length_t(\mu) &= \int_{\mu_s} \frac{\Length(h_{st}(I_a(x)))}{\Length(I_a(x))} \, d\mu_s(x) \\
& =  \int_{\mu_s} S^{t-s}_a(x) \, d\mu_s(x)\\
\intertext{where $S^{r}_a(x)$, the ``stretch'' of the segment $I_a(x)$ with displacement
$r$, is defined to be:}
S^{r}_a(x) &= \frac{\Length(h_{s,s+r}(I_a(x)))}{\Length(I_a(x))} = \int_{\tau \in I_a(x)}
\frac{h_{s,s+r}'(\tau)}{a} \, d\tau
\end{align*}
Now apply this for $r=\pm n$, and we have two versions:
$$\Length_{s+n}(\mu) = \int_{\mu_s} S^{+n}_a(x) \, d\mu_s(x)
$$
and
$$\Length_{s-n}(\mu) = \int_{\mu_s} S^{-n}_a(x) \, d\mu_s(x)
$$
Applying Lemma~\ref{LemmaUniformFlaring}, for each $x \in \mu_s$ the sequence
$\Length(h_{s,s+i}(I_a(x)))$ satisfies $(\kappa,n,A)$ flaring. Taking $a=A$ it follows
that for every $x \in \mu_s$, either $S^{-n}_A(x) \ge \kappa$ or $S^{+n}_A(x) \ge \kappa$.
Now we can define two subsets: 
\begin{align*}
\mu_s^+ = \{x \in \mu_s \suchthat S^{+n}_A(x) \ge \kappa\} \\
\mu_s^- = \{x \in \mu_s \suchthat S^{-n}_A(x) \ge \kappa\}
\end{align*}
Each of these subsets is measurable, and $\mu_s = \mu_s^+ \union \mu_s^-$. It follows that
one of the two sets $\mu_s^+$, $\mu_s^-$ contains at least half of the total $d\mu_s$ measure
(maybe they both do). Choose $\epsilon \in \{+,-\}$ so that $\mu_s^\epsilon$ has more than
half of the measure.

By increasing $n$ to $\left\lceil n \log_\kappa 2 \right\rceil$ if necessary, we may assume
that $\kappa > 2$. It follows that:
\begin{align*}
\Length_{s + \epsilon n}(\mu) &\ge \int_{\mu_s^\epsilon} S^{\epsilon n}_a(x) \, d\mu_s(x) \\
&\ge \kappa \int_{\mu_s^\epsilon} d\mu_s \\
&\ge \frac{\kappa}{2} \Length_s(\mu)
\end{align*}
Since $\kappa/2 > 1$, this proves the lemma.
\end{proof}

\subsection{Growth of measured laminations} 
\label{SectionGrowthOfLaminations}

In addition to the objects fixed at the beginning of Section~\ref{SectionFlaringLaminations},
until further notice we shall fix numbers $L \ge 1$,
$\kappa > 1$, $n \in \Z_+$, and $A \ge 0$, depending only on $\B$, $\rho$, $\delta$, such that
the conclusions of Lemmas~\ref{LemmaLengthFlares} and~\ref{LemmaUniformFlaring} both hold. In
particular, for each $\mu \in \MF$ the sequence $\Length_i(\mu)$, parameterized by
$J=I\intersect\Z$, is $L$--lipschitz and satisfies $(\kappa,n,0)$ flaring. 

Consider now any sequence $\ell_i$, $(i \in J)$, satisfying
$(\kappa,n,0)$ flaring. Let $\kappa'=\kappa^{1/n}$, and recall inequalities~\ref{EquationGrowth}
and~\ref{EquationFullGrowth}. Given $i,j=i+n \in J$, either $\ell_j \ge \ell_i$ or $\ell_i
\ge \ell_j$. In the case $\ell_j \ge \ell_i$ it follows from flaring that $\ell_{j+n}
\ge \kappa \ell_j$, and inequalities~\ref{EquationGrowth}
and~\ref{EquationFullGrowth} hold with $j_0=j$ and $\epsilon=+1$. In the case $\ell_i \ge
\ell_j$ it follows that $\ell_{i-n} \ge \kappa \ell_i$, and the same inequalities hold but with
$j_0=i$ and $\epsilon=-1$. It immediately follows that:

\begin{proposition}
\label{PropTroughLocation}
There exist a constant $\omega$ depending only on $\kappa$, $n$, $L$ (and hence depending
only on $\B$, $\rho$, $\delta$) such that if $(\ell_i)_{i \in J}$ is an $L$--lipschitz sequence
exhibiting $(\kappa,n,0)$ flaring then the following hold.
\begin{itemize}
\item If $\ell_i$ has no minimum for $i \in J$ then $J$ is infinite and $\ell_i$ approaches zero
as $i \to -\infinity$ or as $i \to +\infinity$, but not both; in this case we say that $\ell_i$
achieves its minimum at $-\infinity$ or at $+\infinity$, respectively.
\item If $\ell_i$ has a minimum for $i \in J$, then the smallest subinterval of $J$ containing
all minima, called the \emph{trough} of $\ell_i$, has length at most $\omega$.
\item The sequence grows exponentially as it moves away from the minima, in the following sense:
\begin{itemize}
\item If $i_0 \in J$ is to the right of all minima then for $i \in J$, $i \ge i_0$ we have
$$\ell_i \ge L^{-n} \kappa'{}^{i-i_0} \ell_{i_0} 
$$
and moreover if $i=i_0+kn$ for $k \in \Z_+$ then we have $\ell_i \ge \kappa^n \ell_{i_0}$.
\item If $i_0 \in J$ is to the left of all minima then for $i \in J$, $i \le i_0$ we have
$$\ell_i \ge L^{-n} \kappa'{}^{i_0-i} \ell_{i_0} 
$$
and moreover if $i=i_0-kn$ for $k \in \Z_+$ then we have $\ell_i \ge \kappa^n \ell_{i_0}$.
\end{itemize}
\end{itemize}
\end{proposition}

For each $\mu\in\MF$, Proposition~\ref{PropTroughLocation} applies to the sequence
$\Length_i(\mu)$, $i \in J=I\intersect\Z$. If $\Length_i(\mu)$ achieves its minimum at
$\pm\infinity$, then we say that $\mu$ is \emph{realized at $\pm\infinity$}. On the other hand,
if the length sequence $\Length_i(\mu)$ achieves its minimum at a finite value then we say that
$\mu$ is \emph{realized} at that value; $\mu$ may be realized at several values, and the
\emph{trough} of $\mu$ is defined to be the trough of the length sequence $\Length_i(\mu)$. Since
length is a homogeneous function on $\MF$, these concepts apply as well to elements of $\PMF$.

\begin{corollary}
\label{CorollaryFinitelyRealized}
Each measured lamination consisting of a simple closed geod\-esic is realized at a finite value.
\end{corollary}

\begin{proof}
Since $\B$ is compact and $\gamma$ is $\B$--cobounded, there exists $m>0$ depending only on $\B$
such that for each $t \in I$, every simple closed geodesic in $\S_t$ has length $\ge m$. If $c$
is a simple closed geodesic equipped with a transverse measure $r \in \reals$, it follows that
$\Length_i(c) \ge rm$ for all $i$, and so $c$ cannot be realized at $+\infinity$ or at
$-\infinity$. 
\end{proof}

The next result shows that the position of realizability, either $-\infinity$, a finite
set, or $+\infinity$, is a coarsely continuous function of $\P\mu \in \PMF$. It is a
consequence of the fact that as $\lambda$ varies in $\MF$, the length function
$(\Length_i(\mu))_{i \in J}$ varies continuously in the topology of pointwise convergence.

\begin{proposition} 
\label{PropContinuityOfTrough}
There exists a constant $\Delta=\Delta(\kappa,n,L)$ such that the
following holds for each $\P\mu \in \PMF$.
\begin{enumerate}
\item If $\P\mu$ has trough $W \subset \Z$ then there is a
neighborhood $U \subset \PMF$ of $\P\mu$ such that each $\P\mu' \in U$ has a trough
$W'$ with $\diam(W \union W') \le \Delta$.
\item If $\P\mu$ is realized at $+\infinity$ then for each $i_0\in J$ there exists a neighborhood
$U\subset \PMF$ of $\P\mu$ such that if $\P\mu' \in U$ then either $\P\mu'$ is realized at
$+\infinity$ or $\P\mu'$ is finitely realized with trough contained in $[i_0,+\infinity)$.
\item If $\P\mu$ is realized at $-\infinity$ then a similar statement holds.
\end{enumerate} 
\end{proposition}

\begin{proof} By homogeneity of length it suffices to prove the analogous statement for each
$\mu\in\MF$. 

From the argument preceding the statement of Proposition~\ref{PropTroughLocation}, we immediately
have the following:

\begin{lemma} 
\label{LemmaTroughBounds}
There exists a constant $\zeta \in \Z_+$, depending only on $\kappa$, $n$,
$L$ (and so only on $\B$, $\rho$, $\delta$) such that if $\ell_i$ $(i \in J)$ is $L$--lipschitz and
$(\kappa,n,0)$ flaring then 
\begin{itemize}
\item if $j_0,j_0+n \in J$ and if $\ell_{j_0+n} \ge \ell_{j_0}$ then all minima of $\ell_i$ lie to
the left of $j_0+\zeta$. 
\item if $j_0,j_0-n \in J$ and if $\ell_{j_0-n} \ge \ell_{j_0}$ then all minima of $\ell_i$ lie to
the right of $j_0-\zeta$.\qed 
\end{itemize}
\end{lemma}

To prove (1), let $W=[k,l]$ be the trough of $\mu$. If $k-n\in J$ then $\Length_{k-n}(\mu) >
\Length_k(\mu)$, and if $l+n \in J$ then $\Length_{l+n}(\mu) >
\Length_l(\mu)$. By
continuity of
$\Length\from\T\cross\MF\to (0,\infinity)$ we may choose a neighborhood $U\subset \MF$ of $\mu$
so that if $\mu'\in U$ then
$\Length_{k-n}(\mu') > \Length_{k}(\mu')$ and $\Length_{l+n}(\mu') > \Length_{l}(\mu')$. It
then follows that the trough of $\mu'$ is a subset of the interval
$[k-\zeta,l+\zeta]$, so (1) is proved with $\Delta=\omega+2\zeta$.

To prove (2), assuming $\mu$ is realized at $+\infinity$ it follows that
$\Length_{i_0+\zeta}(\mu) \ge \kappa
\Length_{i_0+\zeta+n}(\mu)$, and so we may choose $U$ so that if $\mu' \in U$ we have
$\Length_{i_0+\zeta}(\mu') > \Length_{i_0+\zeta+n}(\mu')$. It follows that all minima of
$\Length_i(\mu')$ lie to the right of $i_0$. 

The proof of (3) is similar.
\end{proof}

\subsection{Construction of ending laminations} 

We now construct laminations which are realized nearly anywhere one desires, in particular
laminations realized at any infinite ends of $J$. Recall that a measured geodesic lamination is
\emph{perfect} if it has no isolated leaves, or equivalently if it has no
closed leaves.

\begin{proposition}
\label{PropEndLamConstruction}
There exists a constant $\eta$ depending only on $\B$, $\rho$, $\delta$ such that the following
holds. For each $k\in J$ there exists $\mu \in \MF$ which is finitely realized
and whose trough $W$ satisfies $\diam(W \union \{k\}) \le \eta$. If $J$ is infinite then for each
infinite end $\pm\infinity$ of $J$ there exists $\mu^\pm \in \MF$ which is realized at
$\pm\infinity$, respectively; moreover any such $\mu^\pm$ is perfect. 
\end{proposition}

A lamination $\mu^\pm$ realized at an infinite end $\pm\infinity$ is called an \emph{ending
lamination} of $\S_\gamma$. Also, for any finite endpoint $k \in J$, we use the term
\emph{endpoint lamination} to refer to a lamination $\mu$ whose length function $\Length_i(\mu)$
has a minimum occuring with distance $\eta$ of the endpoint $k$; an alternate definition would
require the entire trough of $\mu$ to lie within distance $\eta$ of $k$, but this does not work
out as well, as noted in the remark preceding Proposition~\ref{PropCompactness}.

\begin{proof} As in Corollary~\ref{CorollaryFinitelyRealized}, using compactness of $\B$ and
$\B$--coboundedness of $\gamma$, there exists $m>0$ depending only on $\B$ such that for each $t
\in I$, every simple closed geodesic in $\S_t$ has length $\ge m$. There also exists $M > 0$
depending only on the topology of $S$ such that for each hyperbolic structure on $S$ the shortest
geodesic has length $\le M$; this standard fact follows because the area of any hyperbolic
structure is equal to $2 \pi\chi(S)$, and if the shortest geodesic had arbitrarily large length
then it would have an annulus neighborhood with arbitrarily large area, violating the
Gauss--Bonnet theorem.

Given $k \in J$, take a simple closed geodesic $c$ of minimal length in $\S_k$, with the
transverse Dirac measure. By Corollary~\ref{CorollaryFinitelyRealized}, $c$ is finitely realized.
Consider the subsequence $\Length_{k+np}(c)$, and let $p=P$ be the
value where it achieves its minimum. Since $\Length_{k+nP}(c) \ge m$ and $\Length_k(c) \le M$,
from
$(\kappa,n,0)$ flaring it follows that $\abs{P} < \left\lceil
\log_\kappa
\frac{M}{m} \right\rceil$, and so by Lemma~\ref{LemmaTroughBounds} the trough of
$c$ must be located within the interval
$$\left[k-n\left\lceil \log_\kappa\frac{M}{m} \right\rceil - \zeta \  ,
\   k +n\left\lceil\log_\kappa\frac{M}{m} \right\rceil + \zeta\right]
$$
and so we may take $\eta=2n\left\lceil\log_\kappa\frac{M}{m} \right\rceil + 2\zeta$, proving the
first part of the proposition.

Consider now an infinite end of $J$, say, $+\infinity$. For each $i$ choose $\mu_i \in \MF$ to be
finitely realized, with trough $W_i$ satisfying $\diam(W_i \union \{i\}) \le
\eta$. Using compactness of $\PMF$, choose $\mu^+\in \MF$ so that, after passing to a
subsequence, $\P\mu_i\to \P\mu^+$ as $i \to +\infinity$. 

The fact that $\mu^+$ is realized at $+\infinity$ is a consequence of the fact that the sequence
of length functions $(\Length_j(\mu_i))_{j \in J}$ converges pointwise to the length function
$(\Length_j(\mu^+))_{j \in J}$. To be precise, suppose first that $\mu^+$ is finitely
realized with trough $W \subset \Z$. Applying Proposition~\ref{PropContinuityOfTrough} it follows
that there is a neighborhood $U$ of $\P\mu^+$ and a larger interval $W' \subset \Z$ such that if
$\P\mu' \in U$ then $\P\mu'$ is finitely realized with trough contained in $W'$. But $\P \mu_i
\in U$ for sufficiently large $i$, and its trough $W_i$ goes off to $+\infinity$ as $i \to
+\infinity$, a contradiction. Suppose next that $\mu^+$ is realized at $-\infinity$. By
Proposition~\ref{PropContinuityOfTrough} it follows that the trough of $\P \mu_i$ goes to
$-\infinity$ as $i \to +\infinity$, also a contradiction.

The construction of $\mu^-$ realized at $-\infinity$ is similar.

If $\mu^\pm$ were not perfect it would have a closed leaf $c$ with transverse measure $r \in
\reals$, but then it would follow that $\Length_i(\mu^\pm) \ge rm$ for all $i$, contradicting
that $\Length_i(\mu^\pm) \to 0$ as $i \to \pm\infinity$.
\end{proof}

\subsection{Strict decay of ending laminations} 
\label{SectionStrictDecay}

In this section and the next we concentrate on properties of ending laminations associated to
infinite ends of $\S_\gamma$. The technical lemma~\ref{LemmaStrictDecay} proved in this section
is applied to obtain filling properties for ending laminations. As a consequence, at the end of
section~\ref{SectionFilling}, we will describe the construction of the desired
\Teichmuller\ geodesic $g$ in the case where $\gamma$ is infinite.

Let $\mu^\pm\in\MF$ be an ending lamination realized at an infinite end $\pm\infinity$ of~$J$. Let
$\mu^\pm_i$ denote the measured geodesic lamination on $\S_i$ representing $\mu^\pm$. We prove a
strict decay property for $\mu^+$, say: in any leaf of $\Susp(\mu^+)$, any two connection lines
which are sufficiently far apart in that leaf at level $i$ decay exponentially \emph{immediately}
in the positive direction---there is no growth \emph{anywhere} in the lamination $\mu^+$ as one
approaches $+\infinity$, except on uniformly short segments. A similar statement holds for
$\mu^-$, flowing in the negative direction along connection lines. We make this precise as
follows.

Let $h^\pm_{st}$ be the connection maps on $\Susp(\mu^\pm)$.

\begin{lemma}  
\label{LemmaStrictDecay}
If $\ell$ is a leaf segment of $\mu^+_i$ with $\Length \ell \ge A\kappa$ then  
$$\Length(h^+_{i,i+n}(\ell)) \le \frac{1}{\kappa} \Length\ell.
$$
Similarly, if $\ell$ is a leaf segment of $\mu^-_i$ with $\Length\ell \ge A\kappa$ then
$$\Length(h^-_{i,i-n}(\ell)) \le \frac{1}{\kappa} \Length\ell.
$$
\end{lemma}

\begin{proof}
Borrowing notation from Lemma~\ref{LemmaLengthFlares}, given
$x \in \mu^+_i$ let $I_a(x)$ be the leaf segment of $\mu^+_i$ of length $a$ centered on $x$,
and let $S^r_a(x)$ be the stretch of the segment $I_a(x)$ with displacement $r$, that is:
$$S^r_a(x) = \frac{\Length(h^+_{i+r} I_a(x))}{\Length(I_a(x))}
$$
The lemma says that $S^n_a(x) \le \frac{1}{\kappa}$ if $a \ge A\kappa$. 

If there exists $x \in \mu^+_i$ and $a \ge A\kappa$ such that $S^{n}_{a}(x) >
\frac{1}{\kappa}$, then letting $y$ be the midpoint of $h^+_{i,i+n}(I_a(x))$, and
taking $a' = a \cdot S^{n}_a(x) > A$ we have $S^{-n}_{a'}(y) < \kappa$. By changing variables it
therefore suffices to prove that for all $a \ge A$ and all $x \in \mu^+_i$ we have $S^{-n}_a(x)
\ge
\kappa$.

Suppose there exists $x \in \mu^+_i$ and $a \ge A$ such that $S^{-n}_a(x) < \kappa$. 
Applying Lemma~\ref{LemmaUniformFlaring} we conclude that $S^n_a(x) \ge \kappa$. Now
$S^n_a(y)$ is a continuous function of $y \in \mu^+_i$ and it follows that there is a
neighborhood $U
\subset \mu^+_i$ of $x$ such that if $y \in U$ then $S^n_a(y) > 1$. Given
$y \in U$, again applying Lemma~\ref{LemmaUniformFlaring} it follows by induction on
$p$ that $S^{np}_a(y) \ge \kappa^{p-1}$ for all $p \ge 1$. But $U$ has positive measure
$\int_U d\mu^+_i$ in $\mu^+_i$, and so we have
\begin{align*}
\Length_{i+np}(\mu^+) &= \int_{\mu^+_i} S^{np}_a(y) d \mu^+_i(y) \\
 &\ge \int_U S^{np}_a(y) d \mu^+_i(y) \\
 & \ge \kappa^{p-1} \int_U d\mu^+_i \\
 & \to +\infinity \quad\text{as}\quad p \to +\infinity
\end{align*}
contradicting that $\Length_{i+np}(\mu^+) \to 0$ as $p \to +\infinity$. 
\end{proof}

Recall that two points $\mu,\nu \in \MF$ are \emph{topologically equivalent} if they are
represented by measured foliations which have the same underlying nonmeasured foliation.
Equivalently, for any hyperbolic structure on $S$, the straightenings of $\mu,\nu$ have the
same underlying nonmeasured geodesic lamination.

\begin{corollary}
\label{CorollaryNotTopEquiv} If $J$ is bi-infinite and if $\mu^-$, $\mu^+$ are the ending
laminations realized at $-\infinity$, $+\infinity$ respectively, then $\mu^-$ and $\mu^+$ are not
topologically equivalent.
\end{corollary}

\begin{proof}
Suppose they are topologically equivalent, and so on the surface $\S_0$ the laminations
$\mu^-_0$ and $\mu^+_0$ have the same underlying nonmeasured geodesic lamination. Let $\ell$
be any leaf segment of this lamination with $\Length(\ell) > A\kappa$. Applying
Lemma~\ref{LemmaStrictDecay} twice, from $\mu^+_0$ we conclude that $\Length(h_{0,n}(\ell)) <
\frac{1}{\kappa}
\Length(\ell)$, and from $\mu^-_0$ we conclude that $\Length(h_{0,-n}(\ell)) \le
\frac{1}{\kappa} \Length(\ell)$. However, from Lemma~\ref{LemmaUniformFlaring} at least one
of $\Length(h_{0,n}(\ell))$,
$\Length(h_{0,-n}(\ell))$ is $\ge \kappa \Length(\ell)$, a contradiction.
\end{proof}

\subsection{Individual filling of the ending laminations}
\label{SectionFilling}

Recall that $\mu \in \MF$ \emph{fills} the surface $S$ if $\mu$ has nonzero
intersection number with every simple closed curve. Equivalently, for any hyperbolic
structure, the realization of $\mu$ as a measured lamination has simply connected
complementary components. Note that a filling geodesic lamination is necessarily perfect.

\begin{proposition}
\label{PropEndLamsFill}
Any ending lamination $\mu^\pm$ fills $S$.
\end{proposition}

\begin{proof}[Proof of Proposition \ref{PropEndLamsFill}]
Arguing by contradiction, suppose that, say, $\mu^+$ does not fill $S$. Consider the
straightening $\mu^+_0$ in $\S_0$. Let $F$ be a component of $\S_0 - \mu^+_0$ which is not simply
connected. Let $c$ be a simple closed geodesic which is peripheral in $F$. Let $E$ be a component
of $F-c$ which is a neighborhood of an end of $F$.  The metric completion $\overline E$ is a
``crown'' surface (see Figure~\ref{FigureCrown}), i.e.\ a complete hyperbolic surface with
geodesic boundary homeomorphic to an annulus with $Q \ge 1$ ``crown points'' removed from one of
the boundary components of the annulus. Each removed crown point has a neighborhood isometric to
the region in $\hyp^2$ bounded by two geodesics with a common ideal endpoint in $\bdy\hyp^2$ and
a horocycle attached to that endpoint. The compact boundary component of $\overline E$ is $c$.
The rest of the boundary $\bdy\overline E - c$ consists of $Q$ components, each isometric to the
real line, each identified with a leaf of $\mu^+_0$.

\begin{figure}
\centeredepsfbox{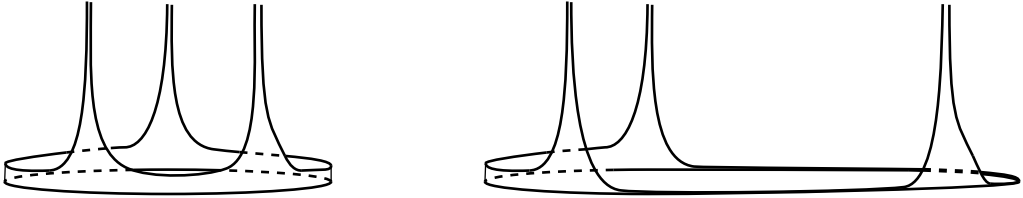}
\caption{If $c$ is a simple closed geodesic in $S$ which is disjoint from a nonfilling geodesic
lamination $\mu$, and if $c$ is peripheral in $S-\mu$, then $S-(\mu\union c)$ has a component $E$
whose metric completion $\overline E$ is a crown. For each $\epsilon>0$, if $c$ is sufficiently
long then the set of points in $c$ that are within distance $\epsilon$ of $\bdy\overline E - c$
consists of at least $1-\epsilon$ of the total length of $c$.}
\label{FigureCrown}
\end{figure}

Let $c_i$ be the straightening of $c$ in $S_i$, let $F_i$ be the component of $S_i -
\mu^+_i$ containing $c_i$, and let $E_i$ be the component of $F_i - c_i$ corresponding to
$E_0$. The metric completion $\overline E_i$ is also a crown surface with $Q$ crown
points, with one compact boundary component $c_i$, and with $\bdy\overline E_i - c_i$ consisting
of $Q$ noncompact boundary components each identified with a leaf of $\mu^+_i$. 

By Corollary~\ref{CorollaryFinitelyRealized}, $c$ is finitely realized, and so $\Length_i(c) \to
\infinity$ as $i \to +\infinity$. It follows that $c_i$ has longer and longer subsegments which
are closer and closer to subsegments of the opposite boundary of $\overline E_i$ (see
Figure~\ref{FigureCrown}). To be precise:

\begin{claim} For each $\epsilon>0$ there exists $i_0$ such that if $i \ge i_0$ then the
set $c^\epsilon_i = \{x \in c_i \suchthat d(x,\bdy \overline E_i-c_i) < \epsilon\}$ consists
of at most $Q$ segments, and $\Length(c^\epsilon_i) \ge (1-\epsilon) \Length(c_i)$. If $\epsilon$
is sufficiently small then each component of $c^\epsilon_i$ is within distance $\epsilon$ of at
most one component of $\bdy\overline E_i - c_i$.
\end{claim}

To see why, the hyperbolic surfaces $\overline E_i$ have constant topology, and therefore they
have constant area, by the Gauss--Bonnet theorem. It follows that $\Length(c_i - c^\epsilon_i)$ is
bounded above, otherwise the $\epsilon$--neighborhood of $c_i - c^\epsilon_i$ would have too much
area. Also, as long as $\epsilon$ is sufficiently small, if $[x,y] \subset c_i$ is a segment such
that $x$ is within $\epsilon$ of one component of $\bdy \overline E_i - c_i$ and $y$ is within
$\epsilon$ of a different component, then there must be a point in $[x,y]$ which has distance $>
\epsilon$ from $\bdy\overline E_i - c_i$, and so $x,y$ are in different components of
$c^\epsilon_i$. In other words, each component of $c^\epsilon_i$ is within distance $\epsilon$ of
only one of the $Q$ components of $\bdy\overline E_i - c_i$.

The idea of the proof of Proposition~\ref{PropEndLamsFill} is that $c_i$ is growing
exponentially, whereas long leaf segments of $\mu_i$ are shrinking exponentially, and since most
of $c_i$ is very close in the tangent line bundle to $\mu_i$ this leads to a contradiction. Now
we make this precise.

As a consequence of the claim, since $\Length(c^\epsilon_i) + \Length(c-c^\epsilon_i) =
\Length(c_i)$, we have
$$\Length(c-c^\epsilon_i) \le \frac{\epsilon}{1-\epsilon} \Length(c^\epsilon_i).
$$
Now choose a very small $\epsilon>0$, and choose $i$ so large that, listing the path components of
$c^\epsilon_i$ as $\alpha_1,\ldots,\alpha_K$, there are corresponding arcs
$\beta_1,\ldots,\beta_K \subset \bdy\overline E_i - c_i \subset
\mu^+_i$, such that for $i=1,\ldots,K$, the arcs $\alpha_k$,
$\beta_k$ are $\epsilon$--fellow travellers, and $\Length(\beta_k) \ge A
\kappa$. Let $\alpha'_k = H_{i,i+n}(\alpha_k) \subset c_{i+n}$, $\beta'_k =
H_{i,i+n}(\beta_k) \subset \mu^+_{i+n}$ where $H_{i,i+n}$ is the connection map on geodesic
laminations, or equivalently, the connection on the geodesic foliation bundle. Applying
Lemma~\ref{LemmaStrictDecay}, we have
$$\Length(\beta_k) \ge \kappa \Length(\beta'_k).
$$
If $\epsilon$ is sufficiently small, each of the pairs of segments $\alpha_k$, $\beta_k$ is
arbitrarily close when lifted to the tangent line bundle, and so by the uniform continuity
property of $H_{i,i+n}$ (see Lemma~\ref{LemmaGFConnection}) they stretch by very nearly equal
amounts:
$$\frac{\Length(\alpha'_k)}{\Length(\alpha_k)} \biggm/ \frac{\Length(\beta'_k)}{\Length(\beta_k)}
\le 1+\eta
$$
for any given $\eta>0$. We therefore have:
\begin{align*}
\Length(c_{i+n}) &= \Length(H_{i,i+n}(c^\epsilon_i)) + \Length(H_{i,i+n}(c-c^\epsilon_i)) \\
                 &\le \Length(H_{i,i+n}(c^\epsilon_i)) + L^n \Length(c-c^\epsilon_i) \\
                 &\le \Length(H_{i,i+n}(c^\epsilon_i)) + L^n \frac{\epsilon}{1-\epsilon}
\Length(c^\epsilon_i) \\
                 &\le \left(1 + \frac{\epsilon L^{2n}}{1-\epsilon} \right)
\Length(H_{i,i+n}(c^\epsilon_i)) \\
                 &= \left(1 + \frac{\epsilon L^{2n}}{1-\epsilon} \right) \Sum_1^K
\Length(\alpha'_k) \\
                 &\le \left(1 + \frac{\epsilon L^{2n}}{1-\epsilon} \right) (1+\eta)  \Sum_1^K
\Length(\alpha_k) \cdot
\frac{\Length(\beta'_k)}{\Length(\beta_k)} \\
                 &\le \left(1 + \frac{\epsilon L^{2n}}{1-\epsilon} \right) \frac{1+\eta}{\kappa} 
\Sum_1^K \Length(\alpha_k) \\
                 &\le \left(1 + \frac{\epsilon L^{2n}}{1-\epsilon} \right) \frac{1+\eta}{\kappa} 
\Length(c_i) \\
                 &\le \left(1 + \frac{\epsilon L^{2n}}{1-\epsilon} \right)
\frac{1+\eta}{\kappa^2} 
\Length(c_{i+n})
\end{align*}
where the last inequality follows from the fact that $\Length(c_{i+n}) \ge \kappa \Length(c_i)$
for sufficiently large $i$.  When $\eta$ and $\epsilon$ are sufficiently small, the
multiplicative constant is arbitrarily close to $1/\kappa^2$, and we obtain a
contradiction.
\end{proof}

When $\gamma$ is bi-infinite we are now in a position to construct the desired \Teichmuller\
geodesic $g$. For any topologically inequivalent pair of laminations, if at least one of them
fills, then the pair jointly fills. We therefore have:

\begin{corollary}
\label{CorollaryFilling}
If $\gamma$ is bi-infinite then any choice of ending laminations for the two ends
jointly fills $S$.
\qed\end{corollary}
 
In the bi-infinite case, with ending laminations $\lambda_1,\lambda_2$, we can therefore define a
\Teichmuller\ geodesic line $g=\geodesic{\P\lambda_1}{\P\lambda_2}$. 

In the half-infinite case, associated to the infinite end there is an ending lamination
$\lambda_1$ which fills, and associated with the finite end there is an endpoint lamination
$\lambda_2$ which is topologically inequivalent to $\lambda_1$, and hence the pair
$\lambda_1,\lambda_2$ jointly fills. Unfortunately we cannot yet prove, when $\gamma$ is finite,
that any pair of endpoint laminations jointly fills---indeed it is not always true without an
extra condition.

\subsection{A compactness property}
\label{SectionCompactnessProp}

In this section we prove a compactness property for ending laminations and endpoint
laminations associated to cobounded, lipschitz paths in \Teichmuller\ space whose associated
hyperbolic plane bundle is a hyperbolic metric space. This will be used in the following section
in two key ways: to prove the desired filling property for finite paths; and to prove
Theorem~\ref{TheoremStableTeichGeod}.

For the last several sections we have been fixing a particular path $\gamma$, but now we want to
let $\gamma$ vary and investigate convergence of the various pieces of geometric data we have
been studying. 

Recall that we have fixed a compact subset $\B\subset\Mod$ and numbers $\rho \ge 1$, $\delta \ge
0$. We also fix a constant $\eta \ge 0$ satisfying the conclusions of
Proposition~\ref{PropEndLamConstruction}, in particular each endpoint lamination is realized
somewhere within distance $\eta$ of the endpoint.

Let $\Gamma_{\B,\rho,\delta,\eta}$ be the set of all triples $(\gamma,\lambda_-,\lambda_+)$
with the following properties: 
\begin{enumerate}
\item $\gamma\from I \to \T$ is a $\B$--cobounded, $\rho$--lipschitz, $\Z$--piecewise affine path
such that $\H_\gamma$ is $\delta$--hyperbolic,
\item $0 \in I$, and each of $\lambda_\pm\in\MF$ is normalized to have length~1 in the hyperbolic
structure $\gamma(0)$,
\item 
\label{ItemRealization}
The lamination $\lambda^+$ is realized in $\S_\gamma$ near the right end, in the following
sense:
\begin{enumerate}
\item If $\gamma$ is right infinite then $\lambda^+$ is realized at $+\infinity$.
\item 
\label{ItemFiniteRealization}
If $\gamma$ is right finite, with right endpoint $M$, then there exists a minimum of the
length sequence $\Length_i(\lambda^+)$ lying in the interval $[M-\eta,M]$.
\end{enumerate}
The lamination $\lambda^-$ is similarly realized in $\S_\gamma$ near the left end.
\end{enumerate}
We give $\Gamma_{\B,\rho,\delta,\eta}$ the product topology, using the usual topology on $\MF$ for the
second and third coordinates $\lambda_-$, $\lambda_+$, and for the $\gamma$ coordinate we use the
compact--open topology. Since the domain interval $I$ may vary, we apply the compact--open topology
to the unique, continuous extension $\gamma\from\reals\to\T$ which is constant on each component
of $\reals-I$.

\paragraph{Remark} The existential quantifier in item~\ref{ItemFiniteRealization} above is
important. In the following proposition, the proof in case~2 would fall apart if
item~\ref{ItemFiniteRealization} were replaced, say, by the statement that the entire trough of
$\lambda^+$ lies in the interval $[M-\eta,M]$.

\begin{proposition}
\label{PropCompactness}
The action of $\MCG$ on $\Gamma_{\B,\rho,\delta,\eta}$ is cocompact.
\end{proposition}

\begin{proof}
Choose a compact subset $\A\subset\T$ such that each $(\gamma,\lambda_-,\lambda_+) \in
\Gamma_{\B,\rho,\delta,\eta}$ may be translated by the action of $\MCG$ so that 
\begin{itemize}
\item[(4)] $\gamma(0)\in\A$.
\end{itemize}
It suffices to prove that the set of $(\gamma,\lambda_-,\lambda_+)$ satisfying (1), (2), (3), and
(4) is compact.

By the Ascoli--Arzela theorem, the set of $\rho$--lipschitz, $\Z$--piecewise affine paths $\gamma
\from I \to \T$ with $\gamma(0)\in\A$ is compact (this is where we use $\Z$--piecewise affine).
The subset of those which are $\B$--cobounded is a closed subset, since $\B$ is closed. The subset
of those for which $\H_\gamma$ is $\delta$--hyperbolic is closed, because if $\gamma_i$ converges
to $\gamma$ then $\H_{\gamma_i}$ converges to $\H_\gamma$ in the Gromov--Hausdorff topology, and
for fixed $\delta$ the property of $\delta$--hyperbolicity is closed in the Gromov--Hausdorff
topology \cite{Gromov:hyperbolic}. 

So far we have we have shown that the set of triples satisfying (1) and (4) is compact, and since
the length function $\T\cross\MF\to(0,\infinity)$ is continuous is follows that the set
satisfying (1), (2), and (4) is compact. It remains to show that the subset of those satisfying
(3) is closed. Let $(\gamma_i,\lambda^-_i,\lambda^+_i)$ be a sequence satisfying (1--4) and
converging to a limit $(\gamma,\lambda^-,\lambda^+)$, necessarily satisfying (1), (2), and (4).
Let $I_i$ be the domain of $\gamma_i$, and $I$ the domain of $\gamma$. We must verify (3), and we
focus on the proof for $\lambda^+$, which will be a consequence of the
continuity of length functions. The detailed proof is broken into cases depending on the nature of
the positive ends of the domain intervals.

\paragraph{Case 1: $I$ is positive infinite} We must prove that $\lambda^+$ is realized at
$+\infinity$ in $\S_\gamma$. If not, then it is realized finitely or at $-\infinity$; pick $i_0
\in I$ so that $\lambda^+$ is realized to the left of $i_0$, either at $-\infinity$ or with
trough to the left of $i_0$. By Proposition~\ref{PropContinuityOfTrough} there is a $\Delta>0$,
depending only on $\B$, $\rho$, $\delta$, such that if $i$ is sufficiently large then
$\lambda^+_i$ is realized to the left of $i_0+\Delta$. If there exist arbitrarily large $i$ for
which $I_i$ is positive infinite then $\lambda^+_i$ is realized at $+\infinity$, an
immediate contradiction. On the other hand, if $I_i$ is positive finite for all sufficiently
large $i$, with right endpoint $M_i$, then the endpoint lamination $\lambda^+_i$ is realized to
the right of $M_i-\eta$, but $M_i$ diverges to $+\infinity$ and so eventually $\lambda^+_i$ is not
realized to the left of $i_0+\Delta$, also a contradiction. 

\paragraph{Case 2: $I$ has finite right endpoint $M$} It follows that for sufficiently large
$i$, the interval $I_i$ also has finite right endpoint $M$, and so each lamination $\lambda_i$ is
realized at some point in the interval $[M-\eta,M]$. By continuity of the length function
$\T\cross\MF\to(0,\infinity)$ it follows that $\lambda^+$ is also realized at some point in this
interval.

\end{proof}

\subsection{Proof of Theorem \ref{TheoremStableTeichGeod}}
\label{SectionProof}

At the end of section~\ref{SectionFilling} we used results about filling to construct the desired
\Teichmuller\ geodesic $g$ in the case where the domain $I$ of the path $\gamma$ is a line. In
the case where $I$ is a segment we need the following:

\begin{proposition}
\label{PropFiniteFilling}
There exists a constant $\Lambda$, depending only on $\B$, $\rho$, $\delta$, $\eta$ such that if
$I=[m,n]$ is a finite segment with $n-m \ge \Lambda$, and if $(\gamma,\lambda^-,\lambda^+) \in
\Gamma_{\B,\rho,\delta,\eta}$ with $\gamma \from I \to \T$, then the $\lambda^-, \lambda^+$
jointly fills $S$.
\end{proposition}

\begin{proof} If there exists no such constant $\Lambda$, then there is a sequence of
examples $(\gamma_i,\lambda^-_i,\lambda^+_i)\in\Gamma_{\B,\rho,\delta,\eta}$ with $\gamma_i \from
I_i \to \T$, such that $\Length(I_i) \to \infinity$, and the pair $\lambda^-_i, \lambda^+_i$ does
not jointly fill. After translating the parameter interval $I_i$ we may assume that $0$ lies
within distance $1/2$ of the midpoint of $I_i$. After acting appropriately by elements of $\MCG$,
we may assume that the sequence $(\gamma_i,\lambda^-_i,\lambda^+_i)$ converges to
$(\gamma,\lambda^-,\lambda^+) \in \Gamma_{\B,\rho,\delta,\eta}$, and it follows that $\gamma$ has
domain $\reals$. By Corollary~\ref{CorollaryFilling}, the pair $\lambda^-,\lambda^+$ jointly
fills. However, the set of jointly filling pairs in $\MF\cross\MF$ is an open subset $\FP$, and so
for sufficiently large $i$ the pair $\lambda^-_i,\lambda^+_i$ jointly fills, a contradiction.
\end{proof}

Now we turn to the proof of Theorem \ref{TheoremStableTeichGeod}.

Consider a $\B$--cobounded, $\rho$--lipschitz, $\Z$--piecewise affine path $\gamma \from I \to
\T$ such that $\H_\gamma$ is $\delta$--hyperbolic. By translating the interval $I$ we may assume
$0 \in I$. Choose $\eta$ satisfying Proposition~\ref{PropEndLamConstruction}, and it follows that
there are $\lambda_-,\lambda_+ \in \MF$ such that $(\gamma,\lambda_-,\lambda_+) \in
\Gamma_{\B,\rho,\delta,\eta}$. Recall that $\lambda^-,\lambda^+$ are normalized to have length~1
in the hyperbolic structure $\gamma(0)$.

Fix a constant $\Lambda$ so that Proposition~\ref{PropFiniteFilling} is satisfied.

First we knock off the case where $\gamma\from I\to\T$ satisfies $\Length(I) < \Lambda$. In this
case $\Length(\gamma) < \rho\Lambda$. Let $g$ be the geodesic segment with the same endpoints as
$\gamma$, and so $\Length(g) \le \Length(\gamma) < \rho\Lambda$. It follows that the Hausdorff
distance between $\image(\gamma)$ and $g$ is at most $\rho\Lambda$. Also, any $\rho$--lipschitz
segment of length $\le\rho\Lambda$ is a $(1,C)$--quasigeodesic with $C=\max\{1,\rho\Lambda\}$.

We may henceforth assume that $\Length(I) \ge \Lambda$.

Now we define the geodesic $g$. For each infinite end of $I$ we have associated an ending
lamination, which determines the corresponding infinite end of $g$; in particular, in the case
where $I$ is bi-infinite we have already defined $g = \geodesic{\P\lambda^-}{\P\lambda^-}$. In
the case where $I$ is half-infinite or finite, we also have a
jointly filling pair $\lambda^-,\lambda^+$ and so we have a geodesic line
$g^*=\geodesic{\P\lambda^-}{\P\lambda^+}$. We must specify a ray or segment on $g^*$, and even in
the bi-infinite case we must specify how the path $\gamma$ is synchronized with this ray or
segment. These tasks are accomplished as follows. 

Recall the notation $\sigma(\lambda,\lambda')$ and $q(\lambda,\lambda')$ for the conformal
structure and quad\-ratic differential determined by a jointly filling pair $\lambda,\lambda' \in
\MF$. For each $t \in I$ we define:
$$a^-(t) = \frac{1}{\Length_t \lambda^-}, \quad a^+(t) = \frac{1}{\Length_t \lambda^+}
$$
These are continuous functions of $t \in I$, and it follows that we have a continuous function
$\Sigma \from I \to g^*$ defined as follows:
$$\Sigma(t)=\sigma(a^-(t)\lambda^-,a^+(t)\lambda^+)
$$
The image of this map is therefore a connected subset of $g^*$ whose closure is the
desired geodesic $g$ (we will in fact show that $\image(\Sigma)$ is closed). We also have a
continuous family of quadratic differentials $Q \from I \to \QD$ defined by
$$Q(t) = q(a^-(t)\lambda^-,a^+(t)\lambda^+)
$$
where $Q(t)$ is a quadratic differential on the Riemann surface $\Sigma(t)$.

Next we prove that the \Teichmuller\ distance between $\gamma(t)$ and $\Sigma(t)$ is bounded
above, by a constant depending only on $\B$, $\rho$, $\delta$ (and $\eta$, which depends in turn
on $\B$, $\rho$, $\delta$). For integer values $t=i$ this follows from the compactness result,
Proposition~\ref{PropCompactness}. To see why, defining $\gamma'(s)=\gamma(s+i)$, the ordered
triple $(\gamma',a^-(i)\lambda^-,a^+(i)\lambda^+)$ lies in the $\MCG$--cocompact set
$\Gamma_{\B,\rho,\delta,\eta}$. The map taking $(\gamma',\lambda'{}^-,\lambda'{}^+) \in
\Gamma_{\B,\rho,\delta,\eta}$ to $(\gamma'(0),\sigma(\lambda'{}^-,\lambda'{}^+)) \in \T\cross\T$
is continuous and $\MCG$--equivariant, and therefore has $\MCG$--cocompact image, and hence the
distance function is bounded above as required. If $t$ is not an integer, there exists an integer
$i$ such that $\abs{t-i}\le 1$, and recalling the lipschitz constant $L$ for length functions $t
\mapsto\Length_t(\lambda)$ it follows that 
$$\abs{\log(a^-(i)/a^-(t))},\abs{\log(a^+(i)/a^+(t))} \le
\log(L)
$$ 
and so $\Sigma(t)$ and $\Sigma(i)$ have \Teichmuller\ distance bounded solely in terms
of $L$ (which depends only on $\B$, $\rho$). Also, $\gamma(t)$ and $\gamma(i)$ have \Teichmuller\
distance at most $\rho$.

Our final task is to prove that $\gamma\from I\to\T$ is a quasigeodesic. Since
$d(\gamma(t),\Sigma(t))$ is bounded in terms of $\B$, $\rho$, $\delta$, it suffices to prove that
the map $\Sigma \from I \to \T$ is a quasigeodesic, with constants depending only on $\B$,
$\rho$, $\delta$. Of course the image of $\Sigma$ is contained in the \Teichmuller\ geodesic
$g^*$, but $\Sigma$ does not have the geodesic parameterization, which it would have had if we
had taken $a^-(t) = e^{-t}$, $a^+(t) = e^t$. Instead, the geodesic parameterization was
sacrificed, and $a^-(t)$, $a^+(t)$ were chosen to guarantee synchronization of $\Sigma$ and
$\gamma$, that is, so that $d(\Sigma(t),\gamma(t))$ is bounded independent of $t$. So, even though
$\Sigma$ is not geodesically parameterized, we can nevertheless show that $\Sigma$ is a
quasigeodesic.

Using the fact that $\image(\Sigma)$ is contained in the geodesic $g$, we can obtain an exact
formula for $d(\Sigma(s),\Sigma(t))$, as follows. Note that $Q(t)$ is not necessarily normalized
so that $\norm{Q(t)}=1$, but we have:
\begin{align*}
\norm{Q(t)} &= \norm{q(a^-(t) \lambda^-,a^+(t) \lambda^+)} \\
            &= a^-(t) a^+(t) \norm{q(\lambda^-,\lambda^+)} \\
            &= a^-(t) a^+(t) \norm{Q(0)}
\end{align*}
The ordered pair of measured laminations $a^-(t)\lambda^-$, $a^+(t)\lambda^+$ can be normalized by
dividing each of them by $\sqrt{\norm{Q(t)}}$, which does not affect $\Sigma(t)$:
\begin{align*}
\Sigma(t) &= \sigma\left(\frac{a^-(t)}{\sqrt{\norm{Q(t)}}} \, \lambda^-,
                                  \frac{a^+(t)}{\sqrt{\norm{Q(t)}}} \, \lambda^+  \right) \\
          &= \sigma\left(\sqrt{\frac{a^-(t)}{a^+(t)}} \, \lambda^-, 
                                             \sqrt{\frac{a^+(t)}{a^-(t)}} \lambda^+ \right)
\end{align*}
It follows that for $s,t \in I$ we have:
\begin{align*}
d(\Sigma(s),\Sigma(t)) &= \frac{1}{2} \abs{\log\left(\frac{a^-(s)}{a^+(s)} \biggm/
                                                 \frac{a^-(t)}{a^+(t)}   \right)} \\
                       &= \frac{1}{2} \abs{\log\left(\frac{a^-(s)}{a^-(t)}\right) 
                                          +\log\left(\frac{a^+(t)}{a^+(s)}\right)}
\end{align*}
Assuming as we may that $s\le t$, we apply Lemma~\ref{LemmaLengthFlares} to obtain a constant $L$
depending only on $\B$, $\rho$, $\delta$ so that 
\begin{align*}
\Length_s(\lambda^+) &\le L^{t-s} \Length_t(\lambda^+) \\
a^+(t) &\le L^{t-s} a^+(s) \\
\log\left( \frac{a^+(t)}{a^+(s)} \right) &\le \log(L) (t-s)
\intertext{and similarly}
\log\left( \frac{a^-(s)}{a^-(t)} \right) &\le \log(L) (t-s)
\intertext{and so}
d(\Sigma(s),\Sigma(t)) &\le \log(L) \abs{t-s}.
\end{align*}
For the lower bound we apply Lemma~\ref{LemmaLengthFlares} again to obtain $(\kappa,n,0)$ flaring
of the sequences $a^+(i)$, $a^-(i)$, with $\kappa$, $n$ depending only on $\B$, $\rho$,
$\delta$. We also use the fact that $a^+(i)$ achieves its minimum near the right end of $I$ and
that $a^-(i)$ achieves its minimum near the left end. To simplify matters, truncate $I$ so that
any finite endpoint of $I$ is an integer divisible by $n$, and if it is a left (resp.\ right)
endpoint $np$ then the minimum of the sequence $(a^-(ni))_i$ (resp.\ $(a^+(ni))_i$) is achieved
uniquely with $i=p$. By Proposition~\ref{PropEndLamConstruction} we need only chop off an amount
of length~$\le\eta+2n$ to achieve this effect; at worst the additive quasigeodesic
constant for $\Sigma$ is increased by an amount depending only on $\B$, $\rho$, $\delta$, and the
multiplicative constant is unchanged. For $s=np < t=nq \in I$ we therefore have
\begin{align*}
\Length_t(\lambda^-) &\ge \kappa^{q-p} \Length_s(\lambda^-) \\
a^-(s)   &\ge \kappa'{}^{q-p} a^+(s), \qquad \kappa'=\kappa^{1/n} \\
\log\left(  \frac{a^-(s)}{a^-(t)} \right)  &\ge \kappa'{}^{t-s} \\
\intertext{and similarly}
\log\left( \frac{a^+(t)}{a^+(s)} \right)  &\ge \kappa'{}^{t-s} \\
\intertext{and so}
d(\Sigma(s),\Sigma(t)) &\ge \log(\kappa') \abs{t-s}.
\end{align*}
For general $s \le t \in I$ pick $s'=np \le t'=nq \in I$ so that $\abs{s-s'}, \abs{t-t'} \le n$,
and we have:
\begin{align*}
d(\Sigma(s),\Sigma(t)) &\ge d(\Sigma(s'),\Sigma(t')) - d(\Sigma(s),\Sigma(s')) -
d(\Sigma(t'),\Sigma(s')) \\
   &\ge \log(\kappa') \abs{t'-s'} - 2n\log(L) \\
   &\ge \log(\kappa') \abs{t-s} - 2n\log(L) - 2n\log(\kappa')
\end{align*}
This completes the proof of Theorem \ref{TheoremStableTeichGeod}.

\section{Model geometries for geometrically infinite ends}
\label{SectionELC}

Throughout this section a \emph{hyperbolic 3--manifold} will always be complete, with finitely
generated, freely indecomposable fundamental group, and with no para\-bolics.

If $N$ is a hyperbolic 3--manifold, Scott's core theorem \cite{Scott:core} produces a \emph{compact
core} $K_N \subset N$, a compact, codimension--0 submanifold whose inclusion $K_N \inject N$ is a
homotopy equivalence. Bonahon proves \cite{Bonahon:ends} that $N$ is \emph{geometrically tame},
which by a result of Thurston \cite{Thurston:GeomTop} implies that the inclusion
$\interior(K_N)\inject N$ is homotopic to a homeomorphism. For each end $e$, let $N_e$ be the
closure of the component of $N-K_N$ corresponding to $e$, and it follows that the inclusion of
the closed surface $S_e = K_N\intersect N_e$ into $N_e$ extends to a homeomorphism $S_e \cross
[0,\infinity) \to N_e$.

The manifold $N$ is the quotient of a free and properly discontinuous action of $\pi_1 N$ on
$\hyp^3$, with limit set $\Lambda \subset S^2=\bdy\hyp^3$ and domain of discontinuity
$D=S^2-\Lambda$. The quotient of the convex hull of $\Lambda$ in $\hyp^3$ is called the convex
hull of $N$, $\Hull(N)=\Hull(\Lambda) / \pi_1 N \subset N$. Let $\overline N = (\hyp^3 \union D) /
\pi_1 N$. 

An end $e$ is \emph{geometrically finite} if $N_e$ is precompact in $\overline N$, in which case
one may choose the homeomorphism $S_e \cross [0,\infinity) \to N_e$ so that it extends to a
homeomorphism $S_e \cross [0,\infinity] \to \overline N_e \subset \overline N$, with $S_e \cross
\infinity$ a component of $\bdy\overline N$. Under the action
of $\pi_1 N$ there is an orbit of components of $D$, each an open disc $\delta$, such that the
projection map $\delta \to S_e \cross \infinity$ is a universal covering map. Since the stabilizer
of $\delta$ in $\pi_1 N$ acts conformally on $\delta$, we obtain a conformal structure on
$S_e \cross \infinity$. This conformal structure is independent of the choice of $\delta$ in the
orbit, and its isotopy class is independent of the choice of the homeomorphism $S_e \cross
[0,\infinity]\to \overline N_e$ extending the inclusion map $S_e\inject N_e$. We therefore obtain
a well-defined point of the \Teichmuller\ space of $S_e$, called the \emph{conformal structure at
$\infinity$} associated to the end $e$.

An end $e$ is \emph{geometrically infinite} if it is not geometrically finite. Bonahon proves
\cite{Bonahon:ends} that if $e$ is geometrically infinite then $e$ is \emph{geometrically tame},
which means that there is a sequence of hyperbolic structures $\sigma_i$ on $S_e$ and pleated
surfaces $g_i \from (S_e,\sigma_i) \to N_e$, each homotopic to the inclusion $S_e \inject N_e$,
such that the sequence of sets $g_i(S_e)$ leaves every compact subset of $N_e$. In this
situation, the $\sigma_i$ form a sequence in the \Teichmuller\ space $\T(S_e)$ which accumulates
on the boundary $\PMF(S_e)$. The unmeasured lamination on $S_e$ corresponding to any accumulation
point gives a unique point in the space $\GL(S_e)$, called the \emph{ending lamination} of $e$.

If $N'$ is another \nb{3}hyperbolic manifold and $f \from N \to N'$ is a homeomorphism, then $f$
induces a bijection between the ends of $N$ and of $N'$, an isometry between the corresponding
\Teichmuller\ spaces, and a homeomorphism between the corresponding geodesic lamination spaces.
Thurston's ending lamination conjecture says that if the end invariants of $N$ and $N'$ agree
under this correspondence, then $f$ is homotopic to an isometry. Define the \emph{injectivity
radius} of a hyperbolic manifold $N$ at a point $x$, denoted $\inj_x(N)$, to be the smallest
$\epsilon\ge 0$ such that the ball of radius $\epsilon$ about $x$ is isometric to the ball of
radius $\epsilon$ in hyperbolic space. The injectivity radius of $N$ itself is $\ds \inj(N) =
\inf_{x\in N} \inj_x N$. The manifold $N$ has \emph{bounded geometry} if $\inj(N) > 0$.

\begin{theorem}[Minsky \cite{Minsky:endinglaminations}]
\label{TheoremMinskyELC}
The ending lamination conjecture holds for complete hyperbolic \nb{3}manifolds $N,N'$ with
finitely generated, freely indecomposable fundamental group and with bounded geometry.
\end{theorem}

In the introduction we recounted briefly how Minsky reduced this theorem to the construction of
model manifolds for geometrically infinite ends, Theorems~\ref{TheoremELCDD}
and~\ref{TheoremELC}. The proofs of these theorems occupy
sections~\ref{SectionPleated}--\ref{SectionSinglyDegenerate}. The construction of a model
manifold for an end of bounded geometry, even when the ambient manifold does not have bounded
geometry, is given in Section~\ref{SectionEndsBoundedGeom}.

\subsection{Pleated surfaces} 
\label{SectionPleated}

We review here several facts about pleated surfaces, their geometry, and their homotopies. For
fuller coverage the reader is referred to \cite{Thurston:GeomTop},
\cite{CanaryEpsteinGreen}, and in particular results of \cite{Minsky:endinglaminations} which are
crucial to our proofs of Theorems~\ref{TheoremELCDD} and~\ref{TheoremELC}.

Given a hyperbolic 3--manifold $N$, a \emph{pleated surface} in $N$ is a $\pi_1$--injective,
continuous map $\theta\from F\to N$ where $F$ is a closed surface, together with a hyperbolic
structure $\sigma$ on $F$ and a geodesic lamination $\lambda$ on $\sigma$, such that
any rectifiable path in $F$ is taken to a rectifiable path in $N$ of the same length, $\theta$ is
totally geodesic on each leaf of $\lambda$, and $\theta$ is totally geodesic on each component of
$F-\lambda$. The minimal such lamination $\lambda$ is called the \emph{pleating locus}
of~$\theta$. We incorporate the hyperbolic structure $\sigma$ into the notation by writing
$\theta\from(F,\sigma)\to N$. 

When $e$ is a geometrically infinite end with neighborhood $N_e \homeo S_e \cross [0,\infinity)$,
we shall assume implicitly that any pleated surface with image in $N_e$ is of the form
$\theta \from (S_e,\sigma) \to N_e$ where $\theta$ is homotopic to the inclusion
$S_e\inject N_e$.  We may therefore drop $S_e$ from the notation and write $\theta \from \sigma
\to N_e$. A similar notational convention will be used when $N \homeo S \cross
(-\infinity,+\infinity)$.

An end $e$ is \emph{geometrically tame} if for each compact subset of $N_e$ there is a
pleated surface $\theta \from \sigma \to N_e$ which misses that compact set. Bonahon proved
\cite{Bonahon:ends} (using free indecomposability of $\pi_1 N$) that each geometrically
infinite end is geometrically tame.

The first lemma controls the geometry of a pleated surface in the large, as long
as the ambient manifold has bounded geometry. Note that any pleated surface is distance
nonincreasing, that is, $(1,0)$--coarse lipschitz.

\begin{lemma} 
\label{LemmaGeomBound}
For each $\epsilon>0$, $g \ge 2$ there exists $d>0$, and a properness gauge
$\rho \from [0,\infinity) \to [0,\infinity)$, such that if $\theta\from (F,\sigma) \to N$ is a
pleated surface with $\genus(F) \le g$ and $\inj(N) \ge \epsilon$, then
\begin{enumerate}
\item \label{ItemDiamBound}
The diameter of $\theta(F)$ in $N$ is $\le d$.
\item \label{ItemUnifProp}
{\rm (\cite{Minsky:endinglaminations}, Lemma 4.4)} The map $\wt\theta \from \wt F \to
\wt N = \hyp^3$ is $\rho$ uniformly proper.
\end{enumerate}
\end{lemma}

The next fact relates the geometry of a nearby pair of pleated surfaces:

\begin{lemma}[\cite{Minsky:endinglaminations}, Lemma 4.5]
\label{LemmaPleatedDistance}
For each $\epsilon>0$, $g \ge 2$, and $a \ge 0$ there exists $r \ge 0$ such that
the following holds. Let $\theta_i \from (F,\sigma_i)\to N$, $i=0,1$, be homotopic pleated
surfaces, and suppose that $\inj(N) \ge \epsilon$, $\genus(F) \le g$, and
$d(\image(\theta_0),\image(\theta_1)) \le a$. Then the distance from $\sigma_0$ to $\sigma_1$ in
the \Teichmuller\ space of $F$ is at most $r$. 
\end{lemma}

Next we need a result controlling the geometry of a homotopy between nearby pleated surfaces.
Recall that if $\T$ is the \Teichmuller\ space of a surface $F$ of genus $g$, then for each
$\epsilon>0$, $r \ge 0$ there exists $\kappa\ge 1$ depending only on $g,\epsilon,r$ such
that if $\sigma_0,\sigma_1 \in \T$ both have injectivity radius $\ge\epsilon$, and if $\sigma_0$,
$\sigma_1$ have distance $\le r$ in $\T$, then there is a $\kappa$--bilipschitz map $\phi \from
\sigma_0 \to\sigma_1$ isotopic to the identity.

\begin{lemma}[\cite{Minsky:endinglaminations}, Lemma 4.2]
\label{LemmaShortHomotopy}
For each $\epsilon>0$, $g \ge 2$, $\kappa \ge 1$ there exists $B \ge 0$ such that the
following hold. If $N$ has injectivity radius $\ge\epsilon$ and $F$ has genus $\le g$, if
$\theta_i \from (F,\sigma_i) \to N$ are homotopic pleated surfaces, and if $\phi \from \sigma_0
\to \sigma_1$ is a $\kappa$--bilipschitz map isotopic to the identity,
then there is a straight line homotopy between $\theta_0$ and $\theta_1\composed\phi$ whose
tracks have length~$\le B$.
\end{lemma}

The next fact gives us the raw material for constructing nearby pleated surfaces as needed. It is
an almost immediate consequence of section~9.5 of \cite{Thurston:GeomTop}; we provide extra
details for convenience.

\begin{lemma} 
\label{LemmaUnifDistr}
There exists a constant $\eta>0$ such that for any homotopy equivalence $S \to N$
from a closed surface $S$ to a complete hyperbolic \nb{3}manifold $N$ without parabolics, every
point in the convex hull of $N$ comes within distance $\eta$ of the image of some pleated surface
$\theta\from S \to N$.
\end{lemma}

\begin{proof} We prove this with a bound $\eta$ equal to the smallest positive real number such
that for any ideal hyperbolic tetrahedron $\tau \subset \hyp^3$ and any two faces
$\sigma_1,\sigma_2$ of $\tau$, any point of $\tau$ is within distance $\eta$ of
$\sigma_1\union\sigma_2$. By slicing $\tau$ with totally geodesic planes passing through the
cusps of $\tau$, it follows that $\eta$ is equal to the thinness constant for the
hyperbolic plane $\hyp^2$, namely $\eta = \log(1+\sqrt{2})$ (as computed in
\cite{Cannon:TheoryHyp}, Theorem 11.8).

For each point $x$ in the convex hull of $N$, there exists a disjoint pair of pleated surfaces
$\theta_0,\theta_1 \from S \to N$ such that if $C$ is the unique component of
$N-(\theta_0(S)\union \theta_1(S))$ whose closure $\overline C$ intersects both $\theta_0(S)$ and
$\theta_1(S)$, then $x \in \overline C$. If the ends of $N$ are $e_0,e_1$, we may choose
$\theta_i$ to separate $x$ from $e_i$, either the convex hull boundary when $e_i$ is geometrically
finite, or a pleated surface very far out when $e_i$ is geometrically tame; we use here 
Bonahon's theorem that $N$, and each end of $N$, is geometrically tame \cite{Bonahon:ends}. It
follows that $x$ is in the image of any homotopy from $\theta_0$ to $\theta_1$.

Now we use the construction of section~9.5 of \cite{Thurston:GeomTop}. Given a homotopic pair of
pleated surfaces $\theta_0,\theta_1 \from S \to N$, this construction produces a homotopy
$\theta_t$, $t\in [0,1]$, such that each map $\theta_t \from S \to N$ in the homotopy is either a
pleated surface, or there exist $t' < t < t''$ in $[0,1]$ such that $\theta_{[t',t'']}$ is an
arbitrarily short straight line homotopy, or $\theta_{[t',t'']}$ is supported on the image of a
locally isometric map $\tau\to N$ where $\tau$ is an ideal tetrahedron in $\hyp^3$. In the latter
situation, there are two cases. In the first case, on page 9.47 of \cite{Thurston:GeomTop}, the
homotopy $\theta_{[t',t'']}$ moves two faces of $\tau$ to the opposite two
faces through the image of $\tau$ in $N$. In the second case, on the bottom of page~9.48 of
\cite{Thurston:GeomTop}, there are actually two tetrahedra involved in the homotopy but we may
homotop through them one at a time; the map $\tau\to N$ identifes two faces of the tetrahedron,
wrapping their common edge infinitely around a closed curve in $N$, and the homotopy
$\theta_{[t',t'']}$ moves one of the remaining two faces to last remaining face through the image
of
$\tau$ in $N$. The upshot in either case is that there are at least two faces of $\tau$ whose
images in $N$ lie on the image of one or the other of the pleated surfaces $\theta_{t'}$,
$\theta_{t''}$, and each point swept out by the homotopy is within a bounded distance of the
union of these two faces. 
\end{proof}

Given a metric space with metric $d$ and two subsets $\alpha,\beta$, let $$d(\alpha,\beta) =
\inf\{d(x,y) \suchthat x\in \alpha, y \in \beta\}.$$ A sequence of homotopic pleated
surface $\theta_i\from\sigma_i \to N_e$, defined for $i \in I$ where $I$ is a subinterval of the
integers, is said to be \emph{uniformly distributed} if there are constants
$B>A>0$ such that for each $i \ge 1$,
$$A \le d(\image(\theta_{i-1}),\image(\theta_i)) \le B
$$
and $\image(\theta_{i})$ separates $\image(\theta_{i-1})$ from $\image(\theta_{i+1})$ in $N$. 

Combining Lemma~\ref{LemmaGeomBound}(\ref{ItemDiamBound}) with Lemma~\ref{LemmaUnifDistr} we
immediately obtain:

\begin{lemma}
\label{LemmaPleatedSequence} Suppose $N$ has bounded geometry.
Each geometrically infinite end $N_e$ has a uniformly distributed sequence of pleated surfaces
$\theta_i \from \sigma_i \to N_e$, $i \ge 1$, escaping the end. If $N \homeo S \cross
(-\infinity,+\infinity)$ has two geometrically infinite ends then there is a uniformly
distributed sequence of pleated surfaces $\theta_i \from \sigma_i \to N_e$, $i \in \Z$, escaping
both ends. 
\end{lemma}

In the first case of this lemma, note that if $\wt N \to N$ is the covering space corresponding to
the injection $\pi_1 S_e \inject \pi_1 N$, then the sequence $\theta_i \from \sigma_i \to N_e$
must be contained in the projection of $\Hull(\wt N)$ to $N$.

\subsection{The doubly degenerate case: Theorem \ref{TheoremELCDD}}
\label{SectionDD}

Let $N$, as above, be a complete hyperbolic \nb{3}manifold with finitely generated, freely
indecomposable fundamental group and with bounded geometry. Let $e$ be a geometrically
infinite end of $N$ with corresponding surface $S=S_e$ and neighborhood $N_e \approx S \cross
[1,+\infinity)$, and we assume that the injection $\pi_1(S) \inject \pi_1(N)$ is doubly
degenerate. Let $\wt N \to N$ be the covering map corresponding to the subgroup $\pi_1(S)$, so we
have a homeomorphism $\wt N \approx S \cross (-\infinity,+\infinity)$. The manifold $\wt N$ has an
end $\tilde e$ with neighborhood $\wt N_{\tilde e}$ such that the covering map restricts to an
isometry $\wt N_{\tilde e} \to N_e$, and we may assume that the notation is chosen so that this
isometry is expressed by the identity map on $S \cross [1,+\infinity)$. 

By Theorem~9.2.2 of \cite{Thurston:GeomTop}, each end of $\wt N$ has a neighborhood mapping
properly to a neighborhood of a geometrically infinite end of $N$, and hence $N$ has at most two
ends, each geometrically infinite. 

Consider first the case that $N$ has exactly two ends, and so the covering map
$\wt N \to N$ is bijective on ends. It follows clearly that the covering map $\wt N\to N$
has degree~1 on some neighborhood of $e$. Since degree is locally constant, the covering map is a
homeomorphism, and we have $N\approx \wt N \approx S \cross (-\infinity,+\infinity)$. 

Consider next the case that $N$ has just the one end $e$. The manifold $\wt N$ has two ends,
$\tilde e$ with neighborhood $\wt N_{\tilde e} \approx S \cross [1,+\infinity)$, and $\tilde e'$
with neighborhood $\wt N_{\tilde e'} \approx S \cross (-\infinity,-1]$. Each of these maps by a
finite degree covering map to $N_e \approx S \cross [1,+\infinity)$, in each case induced by a
covering map $S \to S$ which must be degree~1, and hence $\wt N \to N$ is a degree~2 map. It
follows that $N$ has a compact core which is doubly covered by $S \cross [-1,+1]$, and which
has an orbifold fibration with generic fiber $S$ and base orbifold $[0,1]$ with a $\Z/2$ mirror
point at $0$. This implies that $N$ itself has an orbifold fibration with generic fiber $S$ and
base orbifold $[0,\infinity)$ with a $\Z/2$ mirror point at $0$.

Here for convenience is a restatement of theorem~\ref{TheoremELCDD}:

\begin{theorem*}[Doubly degenerate model manifold] Under the above conditions,
there exists a unique cobounded geodesic line $g$ in $\T(S)$ such that $\Sigma(S \to N)$ is
Hausdorff equivalent to $g$ in $\T(S)$. Moreover: the homeomorphism $S \cross
(-\infinity,+\infinity) \approx \wt N$ is properly homotopic to a map which lifts to a
quasi-isometry $\H^\solv_g\to\hyp^3$; and the ideal endpoints of $g$ in $\PMF(S)$ are the
respective ending laminations of the two ends of~$\wt N$. In the degree~2 case, the order~2
covering transformation group on $\wt N$ acts isometrically on $\S^\solv_g$ so as to commute with
the homeomorphism $\S^\solv_g\approx \wt N$.
\end{theorem*}

We'll focus on the two-ended case $N \approx S \cross (-\infinity,+\infinity)$, mentioning later
the changes needed for the one-ended case.

Let $\T$ be the \Teichmuller\ space of $S$.

Applying Lemma~\ref{LemmaPleatedSequence}, let $\theta_n \from (S,\sigma_n) \to N$, $n \in \Z$, be
a uniformly distributed sequence of pleated surfaces. By Lemma~\ref{LemmaPleatedDistance} there
is a constant $\rho$ such that for any $n \in \Z$, the distance in $\T$ between $\sigma_n$ and
$\sigma_{n+1}$ is at most~$\rho$. It follows that there is a $\Z$--piecewise affine,
$\rho$--lipschitz path $\gamma\from \reals \to \T$ with $\gamma(n)=\sigma_n$. Since
$\inj(\sigma_n)$ is uniformly bounded away from zero it follows that the path $\gamma$ is
cobounded. 

Consider now the canonical hyperbolic surface bundle $\S_\gamma \to \reals$ and its universal
cover, the canonical hyperbolic plane bundle $\H_\gamma \to \reals$. We have an identification 
$\S_n \approx\sigma_n\in\T$ for every $n \in \Z$. 

\begin{claim} 
\label{ClaimQILift}
There exists a map $\Phi\from\S_\gamma\to N$ in the correct proper homotopy class which lifts to a
quasi-isometry of universal covers
$\wt\Phi\from\H_\gamma \to \hyp^3$. 
\end{claim}

\begin{proof} To construct the map $\Phi$, restricting the bundle $\S_\gamma$ to
the affine subpath $\gamma \restrict [n,n+1]$ we obtain the bundle $\S_{[n,n+1]}\to [n,n+1]$, on
which there is a connection with bilipschitz constant $\kappa$ depending only the coboundedness
and lipschitz constant of $\gamma$. In particular, we obtain a $\kappa$--bilipschitz map $\phi_n
\from \sigma_n
\to \sigma_{n+1}$, which may be regarded as a map on $S$ isotopic to the identity; but in fact,
under the topological identification $\sigma_n\approx \S_n \approx S$, the map
$\phi_n \from S \to S$ must be the identity. Applying Lemma~\ref{LemmaShortHomotopy} there is a
constant
$B$ and, for each
$n$, a straight line homotopy in $N$ with tracks of length $\le B$ from $\theta_n$ to
$\theta_{n+1}$. The domain of this homotopy may be taken to be $\S_{[n,n+1]} \approx S \cross
[n,n+1]$, and in particular the homotopy is a $B$--lipschitz map when restricted to any $x
\cross [n,n+1]$. Piecing these homotopies together for each $n$ we obtain the desired map
$\Phi\from\S_\gamma\to N$. 

To prove that any lift $\wt\Phi \from \H_\gamma \to \hyp^3$ is a quasi-isometry, obviously
$\wt\Phi$ is surjective, and so  it suffices by Lemma~\ref{LemmaQIFacts}(\ref{ItemUnifPropIsQI})
to prove that $\wt\Phi$ is coarse lipschitz and uniformly proper. Coarse lipschitz is immediate
from the fact that $\wt\Phi \restrict \H_n$ is distance nonincreasing for each $n$, and that
$\wt\Phi$ is $B$--lipschitz along connection lines (indeed this implies that $\wt\Phi$ is
lipschitz). 

Uniform properness of $\wt\Phi$ will follow from two facts.
Lemma~\ref{LemmaGeomBound}(\ref{ItemUnifProp}) tells us that the maps $\wt\theta_n
\from \wt\sigma_n \to \hyp^3$, which are identified with the maps $\wt\Phi \restrict \H_n$, are
uniformly proper with a properness gauge independent of~$n$. Also, the images $\wt\Phi(\H_n)$ are
uniformly distributed in $\H_\gamma$, by Lemma~\ref{LemmaPleatedSequence}. We put these together
as follows.

Consider a number $r \ge 0$ and points $x,y \in \H_\gamma$. We must show that there is a number
$s$ independent of $x,y$, depending only on $r$, such that
\begin{equation}
\text{if}\quad d(x,y) \ge s \quad\text{then}\quad d(\xi,\eta) \ge r,
\quad\text{with}\quad
\xi = \wt\Phi(x), \eta = \wt\Phi(y)
\label{InequalityProper}
\end{equation}
We claim there are constants $L \ge 1, C \ge 0$ such that if $x \in \H_t$ and $y \in \H_u$, then
$d(\xi,\eta) \ge \frac{1}{L} \abs{t-u} - C$. To see why, note that we may assume that $t,u \in \Z$
and $t<u$. Consider the geodesic $\overline{\xi\eta}$ in $\hyp^3$. By
Lemma~\ref{LemmaPleatedSequence}, along this geodesic there is a monotonic sequence of points
$\xi=\zeta_t,\zeta_{t+1},\ldots,\zeta_{u}=\eta$ such that $\zeta_n\in \wt\Phi(\H_n)$, and such
that $\Length(\overline{\zeta_n \zeta_{n+1}})$ is bounded away from zero, establishing the claim.

In order to prove \ref{InequalityProper} it therefore suffices to consider the case $x \in \H_t$,
$y \in \H_u$ with $\abs{u-t} \le R_0=Lr+LC$. Let $y'$ be the point of $\H_t$ obtained by moving
from $y$ into $\H_t$ along a connection line, and so $d(y,y') \le R_0$. Setting
$\eta'=\wt\Phi(y')$ we have $d(\eta,\eta') \le BR_0$. To prove that $d(\xi,\eta) \ge r$ it
therefore suffices to prove
$$d(\xi,\eta') \ge R_1=r+BR_0.
$$
The restriction of $\wt\Phi$ to $\H_n$ is just a lift of the
pleated surface $\theta_n \from \sigma_n \to N$. Applying
Lemma~\ref{LemmaGeomBound}(\ref{ItemUnifProp}) it follows that $\wt\Phi
\restrict \H_n$ is $\rho$-uniformly proper with $\rho$ independent of $n$, and so
$$d(\xi,\eta') \ge \rho(d_n(x,y'))
$$
where $d_n$ is the distance function on $\H_n$ (an isometric copy of $\hyp^2$). By properness of
$\rho$ there exists
therefore a number $R_2 \ge 0$ depending only on $R_1$ such that
$$\text{if}\quad d_n(x,y') \ge R_2 \quad\text{then}\quad\rho(d_n(x,y')) \ge R_1
$$
and so it suffices to prove $d_n(x,y') \ge R_2$. Since the inclusion map $\H_n \inject \H_\gamma$
is distance nonincreasing it suffices to prove $d(x,y') \ge R_2$. Setting $s=R_2 + R_0$ and using
the fact that $d(y,y') \le R_0$, we have:
$$\text{if}\quad d(x,y) \ge s \quad\text{then}\quad d(x,y') \ge R_2
$$
Putting it altogether, if $d(x,y) \ge s$ then $d(\xi,\eta) \ge r$. This completes the proof that
$\wt\Phi$ is a quasi-isometry.
\end{proof}

Applying Claim~\ref{ClaimQILift}, since $\hyp^3$ is a hyperbolic metric space, it follows that
$\H_\gamma$ is also a hyperbolic metric space. Theorem~\ref{TheoremStableTeichGeod} now applies,
and we conclude that there is a bi-infinite geodesic $g$ in $\T$ such that $\gamma$ and $g$
are asynchronous fellow travellers, that is, there is a quasi-isometry $s \from \reals \to
\reals$ such that $d(\gamma(t),g(s(t)))$ is bounded independent of $t$. Note that since $\gamma$
is cobounded and $\gamma$ and $g$ are asynchronous fellow travellers, it follows that $g$ is
cobounded. Applying Proposition~\ref{PropMetricPerturbation} it follows that the map $t
\mapsto s(t)$ lifts to a map $\S^\solv_g \to \S_\gamma$ any of whose lifts $\H^\solv_g \to
\H_\gamma$ is a quasi-isometry. By composition we therefore obtain a map $\S^\solv_g \to N$ any of
whose lifts to universal covers
$\H^\solv_g\to \hyp^3$ is a quasi-isometry, as required for proving Theorem~\ref{TheoremELCDD}.

To complete the proof in the case where $N \homeo S \cross
(-\infinity,+\infinity)$, there are a few loose ends to clean up.

From what we have proved it is evident that $g$ is contained in a bounded neighborhood of the
subset $\{\sigma_i\}$ in $\T$, and hence $g$ is contained in a bounded neighborhood of the subset
$\Sigma(S\to N)$. The opposite containment is an immediate consequence of the fact that the
pleated surface sequence $\theta_i \from \sigma_i \to N$ is uniformly spaced, combined with
Lemma~\ref{LemmaGeomBound}(\ref{ItemDiamBound}) and Lemma~\ref{LemmaPleatedDistance}.

Since $\gamma$ is evidently cobounded, it follows that $g$ is cobounded.

Uniqueness of $g$ is a consequence of the general fact that if $g,g'$ are two cobounded
\Teichmuller\ geodesics whose Hausdorff distance is finite then $g=g'$; for a proof see
\cite{FarbMosher:quasiconvex} Lemma~2.4. 

The fact that the two ends of the geodesic $g$ are the two ending laminations of $N$ is a
consequence of the following fact: if $\theta_i \from \sigma_i \to N$, $i \ge 1$, is a sequence of
pleated surfaces going out an end $e$ of $N$, then for any sequence $\sigma'_i \in \T(S)$ such
that $d(\sigma_i,\sigma'_i)$ is bounded, any accumulation point of $\sigma'_i$ in
$\overline\T(S)$ is an element of $\PMF(S)$ whose image in $\GL(S)$ is the ending
lamination of $e$. For the proof see Lemmas~9.2 and~9.3 of \cite{Minsky:endinglaminations}.

\bigskip

Now we turn to the case where $N$ fibers over the ray orbifold with generic fiber $S$
and one geometrically infinite end. The above analysis applies to the degree~2 covering
manifold $\wt N \homeo S \cross (-\infinity,+\infinity)$, for which we obtain a cobounded
\Teichmuller\ geodesic $g$ and a map $\S^\solv_g \to \wt N$. There is an order~2 isometric
covering transformation $\tau \from \wt N \to \wt N$ which exchanges the two ends, inducing an
order~2 mapping class on $S$ and an order~2 isometry of the \Teichmuller\ space~$\T$. We may then
easily tailor the proof of Lemma~\ref{LemmaUnifDistr} to produce a uniformly distributed sequence
of pleated surfaces $\theta_i \from \sigma_i \to N$ such that in $\T$ we have
$\tau(\sigma_i)=\sigma_{-i}$ for all $i\in \Z$. By the uniqueness clause of
Theorem~\ref{TheoremELCDD} for $\wt N \homeo S \cross (-\infinity,+\infinity)$ it follows that
$\tau$ preserves the \Teichmuller\ geodesic $g$, acting on it by reflection across some fixed
point. This implies that the homeomorphism $\S^\solv_g \to\wt N$ may be chosen in its proper
homotopy class so that the action of $\tau$ on $\wt N$ commutes with an action on $\S^\solv_g$.

\subsection{The singly degenerate case: Theorem \ref{TheoremELC}} 
\label{SectionSinglyDegenerate}

Again let $N$ be a complete hyperbolic \nb{3}manifold with finitely generated, freely
indecomposable fundamental group and with bounded geometry, and let $e$ an end of $N$
with corresponding surface $S=S_e \inject N$ bounding a neighborhood $S \cross [0,\infinity)
\approx N_e \subset N$ of the end $e$. We assume that the injection $\pi_1(S) \inject \pi_1(N)$
is singly degenerate. By lifting to the covering space of $N$ corresponding to the subgroup
$\pi_1(S) \inject N$ we may assume that we have $N\approx S\cross (-\infinity,+\infinity)$. Since
we are in the singly degenerate case, the convex hull of $N$ is bounded by a surface homotopic to
$S$, and so we may assume that $S$ is equal to this surface, and hence $N_e$ is the convex hull
of $N$.

Next we proceed as in the doubly degenerate case to get a uniformly distributed sequence of
pleated surfaces $\theta_n \from (S,\sigma_n) \to N_e$, $n \ge 0$, where $\theta_0 \from S \to
N_e$ is just the inclusion. We obtain a $\Z$--piecewise affine, cobounded, lipschitz ray $\gamma
\from [0,\infinity) \to \T$ with $\gamma(n)=\sigma_n$. The proof of Claim~\ref{ClaimQILift} goes
through as before, and we obtain a map $\Phi \from \S_\gamma \to N_e$ in the correct proper
homotopy class any of whose lifts $\wt\Phi \from \H_\gamma \to \wt N_e$ is a quasi-isometry,
where $\wt N_e$ is the universal cover of $N_e$.

Now there is a trick: the convex hull of the limit set of $\pi_1 S$ is a convex subset of
$\hyp^3$. This convex hull is precisely the space $\wt N_e$, and hence the restriction of the
geodesic metric on $\hyp^3$ is a hyperbolic geodesic metric on $\wt N_e$. The metric space
$\H_\gamma$ is therefore also hyperbolic, and so we may apply Theorem~\ref{TheoremStableTeichGeod}
to $\gamma$, concluding that there is a geodesic ray $g\from [0,\infinity) \to \T$ such that
$\gamma$ and $g$ are asynchronous fellow travellers, and such that $g(0)=\gamma(0)$. 

The proof that $\Sigma(S \to N)$ is Hausdorff equivalent to $g$ in $\T$ goes through as in the
previous case. Given that $g$ has base point $\gamma(0)$, uniqueness of $g$ is proved as before
using coboundedness of $g$. For any other base point $\tau \in \T$, there is a unique ray
$g_\tau$ in $\T$ with base point $\tau$ asymptotic to $g$; uniqueness of $g_\tau$ is proved as
before, and existence is a standard fact. 

Since $\gamma$ and $g_\tau$ are Hausdorff equivalent and $\gamma$ is a quasigeodesic, it follows
that $\gamma$ and $g_\tau$ are asynchronous fellow travellers. The same argument as before
produces a proper homotopy equivalence $N_e \to \S_{g_\tau}$ in the correct proper homotopy
class which lifts to a quasi-isometry $\wt N_e \to \H_{g_\tau}$. 

The proof that the endpoint of $g$ in $\PMF(S)$ gives the ending lamination of $e$ is as before.

\subsection{Ends of bounded geometry} 
\label{SectionEndsBoundedGeom}
Suppose that $N$ is a hyperbolic 3--manifold as
above, but without bounded geometry. Given an end $e$ of $N$, we say that $e$ has bounded
geometry if there exists $\epsilon>0$ such that $\inj(N_e) \ge \epsilon$, meaning that for each
$x \in N_e$ we have $\inj_x(N) \ge \epsilon$. In this case we can push the construction of the
model manifold through; we are grateful to Jeff Brock for suggesting this possibility. We shall
merely sketch the outline of how to do this.

The first tricky part is that the pleated surface results of Minsky quoted in
Section~\ref{SectionPleated} must be adapted to hold without the assumption that $N$ has bounded
geometry. Instead, using only bounded geometry of the end $e$, one must prove that the results
hold in $N_e$.  

In the statement of Lemma~\ref{LemmaGeomBound}, the hypothesis that $\inj(N)>\epsilon$ is replaced
with $\inj(N_e)>\epsilon$, and the conclusion that $\tilde\theta\from \tilde F \to \hyp^3$ is
$\rho$--uniformly proper is replaced by the conclusion that $\tilde\theta \from \tilde F \to \tilde
N_e$ is $\rho$--uniformly proper. With these changes, Minsky's proof taken from
\cite{Minsky:endinglaminations} Lemma~4.4 (which goes back to \cite{Minsky:GeodesicsAndEnds}
Lemma~4.5)  goes through, using basic compactness results for pleated surfaces such as the
theorem of Thurston quoted in \cite{Minsky:GeodesicsAndEnds} as Theorem~4.1.

The statement of Lemma~\ref{LemmaPleatedDistance} undergoes a similar change, and again Minsky's
proof taken from \cite{Minsky:endinglaminations} Lemma~4.5 (which goes back to
\cite{Minsky:GeodesicsAndEnds} Lemma~4.6) goes through. 

Lemma~\ref{LemmaShortHomotopy} also undergoes a similar change.

With these results suitably adapted, it follows that there is an evenly spaced sequence of
pleated surfaces $\theta_i \from \sigma_i \to N_e$ escaping the end $e$, and there is a
$\Z$--piecewise affine, cobounded, Lipschitz ray $\gamma \from [0,\infinity) \to \T$ with
$\gamma(i)=\sigma_i$. The proof of Claim~\ref{ClaimQILift} goes through, producing a map
$\Phi\from\S_\gamma\to N_e$ in the correct proper homotopy class, which lifts to a quasi-isometry
$\tilde\Phi \from\H_\gamma \to \wt N_e$. 

Now we need to check that $\H_\gamma$, or equivalently $\wt N_e$, is a hyperbolic metric space.
The convexity trick used in Section~\ref{SectionSinglyDegenerate} is not available; instead, we
appeal to the methods of \cite{FarbMosher:quasiconvex}. As remarked earlier, hyperbolicity of
$\H_\gamma$ is equivalent to the horizontal flaring property, with a uniform relation between the
constants. It therefore suffices to show directly that the horizontal flaring property holds in
$\H_\gamma$. Suppose that $\alpha,\beta \from I \to \H_\gamma$, $I \subset [0,\infinity)$,
are $\lambda$--quasihorizontal paths in $\H_\gamma$. We must prove that the
sequence $d_i(\alpha(i),\beta(i))$ satisfies $\kappa,n,A(\lambda)$ flaring with appropriate
constants, and we do this by applying the proof of Lemma~5.2 of \cite{FarbMosher:quasiconvex}.
All that is needed to apply that proof is to show that a $\lambda$--quasihorizontal path in
$\H_\gamma$ maps via $\Phi$ to a quasigeodesic in $\hyp^3$, with quasigeodesic constants
depending only on $\lambda$; but this follows immediately from even spacing of the pleated
surface sequence $\theta_i \from \sigma_i \to N_e$. 

Since $\H_\gamma$ is a hyperbolic metric space, Theorem~\ref{TheoremStableTeichGeod} applies as
above to conclude that $\gamma$ fellow travels a unique cobounded geodesic ray $g$ in $\T$ with
$\gamma(0)=g(0)$, and from Proposition~\ref{PropMetricPerturbation} we obtain the desired
quasi-isometry $\H^\solv_g \to \H_\gamma \xrightarrow{\tilde\Phi} \wt N_e$.

\vglue -30pt
\noindent$\phantom{XXX}$

\end{document}